\documentclass[reqno,makeidx,12pt,a4paper]{amsart}
        
        \newif\ifpdf 
        \ifx\pdfoutput\undefined 
        \pdffalse 
        \else 
        \pdfoutput=1 
        \pdftrue 
        \fi

\usepackage{graphicx}
\usepackage{latexsym}
\usepackage{amstext}
\usepackage {amsmath}
\usepackage {amsfonts}
\usepackage {amssymb}
\usepackage {amsthm}
\ifx\cs\documentclass \footheight 12pt \fi \footskip 30pt 
\usepackage{enumerate}
\usepackage{array}
\usepackage{xspace}

\DeclareMathOperator{\dist}{dist}
\DeclareMathOperator{\Id}{Id}

\DeclareMathOperator{\spt}{supp}
\DeclareMathOperator\supp{supp}
\DeclareMathOperator{\dive}{div}
\def\loc{{\mathrm{loc}}}
\DeclareMathOperator{\Lp}{L}

\newcommand{\eps}{\varepsilon}

\newcommand{\Chi}{\mathcal{X}}
\newcommand{\R}{\ensuremath{\mathbb{R}}}
\newcommand{\Rn}{\ensuremath{{\mathbb{R}^d}}}

\newcommand{\N}{\ensuremath{\mathbb{N}}}

\newcommand{\Ha}{\ensuremath{\mathcal{H}}}

\newcommand{\schwto}{\ensuremath{\rightharpoonup}}

\newcommand{\sd}{\, d}
\newcommand{\M}{\mathcal{M}}
\newcommand{\Div}{\operatorname{div}}
\newcommand{\STC}{\kappa}
\newcommand{\ol}{\overline}
\newcommand{\weight}[1]{\left\langle #1 \right\rangle}

\newcommand{\RM}{\lambda}

\def\step#1{\newline\underline{Step\ #1}:\ }

\def\R{\mathbb R}

\def\Per{\mathcal{P}}

\def\per#1{\int_\Omega d|\nabla #1|}
\newcommand{\dprod}[2][]{\langle #2 \rangle_{#1}}
\newcommand{\bdprod}[2][]{\big\langle #2 \big\rangle_{#1}}

\def\BV{\mathrm{BV}}

\theoremstyle{plain}
\numberwithin{equation}{section}
\newtheorem{lemma}{Lemma}[section]
\newtheorem{theorem}[lemma]{Theorem}
\newtheorem{proposition}[lemma]{Proposition}
\newtheorem{definition}[lemma]{Definition}

\theoremstyle{definition}
\newtheorem{remark}[lemma]{Remark}
\def\Xns{H^1_{0,\sigma}(\Omega)}
\def\Xms{BV_{(m_0)}(\Omega;\{0,1\})}
\def\Xmu{H^1_{(0)}(\Omega)}
\def\Xmustar{\Xmu^*}
\def\cp{\mu}
\def\vf{\vartheta}%
\begin{document}
\title[Non-classical sharp interface model for two-phase flows]{Existence of weak solutions for a
non-classical sharp interface model for a two-phase flow of viscous,
  incompressible fluids}
\author{Helmut Abels}
\address{Helmut Abels,  Max Planck Institute for Mathematics in the
  Sciences, Inselstr. 22, D-04103 Leipzig}

\author{Matthias R{\"o}ger}
\address{Matthias R\"{o}ger, Max Planck Institute for Mathematics in the
  Sciences, Inselstr. 22, D-04103 Leipzig}

\email{abels@mis.mpg.de, roeger@mis.mpg.de}



\date{\today}

\begin{abstract}
We introduce a new sharp interface model for the flow of two
immiscible, viscous, incompressible fluids. In contrast
to classical models for two-phase flows we prescribe an evolution law
for the interfaces that takes diffusional effects into
account. This leads to a coupled system of Navier--Stokes and 
Mullins--Sekerka type parts that coincides with the asymptotic
limit of a diffuse interface model.  We
prove the long-time 
existence of weak solutions, which is an open problem for the classical
two-phase model. We show that the phase interfaces have in almost all points
a generalized mean curvature. 
\end{abstract}

\maketitle
\section{Introduction}

We study the flow of two incompressible, viscous and immiscible fluids like
oil and water inside a bounded domain $\Omega$ or in $\Omega=\mathbb{T}^d$,
$d=2,3$. The fluids fill domains $\Omega^+(t)$ and $\Omega^-(t)$, $t>0$,
respectively, with a common interface $\Gamma(t)$ between both
fluids. The flow is 
described in terms of the velocity $v\colon \Omega\times (0,\infty) \to \R^d$ and
the pressure $p\colon \Omega \times (0,\infty) \to \R$ in both fluids in Eulerian coordinates.
We assume the fluids to be of Newtonian type, i.e., the stress tensors are of
the form $T^\pm (v,p)= \nu^\pm Dv -pI$ with constant viscosities $\nu^\pm>0$ and
$2Dv= \nabla v +\nabla v^T$. Moreover, we consider the case with surface
tension at the interface and assume that the densities are the same (and set
to $1$ for simplicity). For the evolution of the phases we take diffusional
effects into account and consider a contribution to the flux that is
proportional to the negative gradient of the chemical potential $\mu$.
Precise assumptions are made below.

To formulate our model we introduce some notation first. Denote by $n$ the unit
normal of $\Gamma(t)$ that points inside $\Omega^+(t)$ and by $V$ and
$H$ the normal velocity and scalar mean curvature of $\Gamma(t)$
with respect to $n$. By $[\cdot]$ we denote
the jump of a quantity across the interface in direction of $n$, i.e.,
$[f](x)= \lim_{h\to 0}(f(x+hn)-f(x-hn))$ for $x\in \Gamma(t)$. 
Then our model is described by the following equations 
\begin{alignat}{2}\label{eq:1}
  \partial_t v + v\cdot \nabla v - \dive T^\pm(v,p) &= 0 &\quad & \text{in}\
  \Omega^\pm(t), t>0, \\\label{eq:2}
\dive v &= 0 &\ & \text{in}\ \Omega^\pm(t), t>0, \\\label{eq:2'}
m\Delta \mu &= 0 &\ & \text{in}\ \Omega^\pm(t), t>0, \\\label{eq:3}
-n\cdot [T(v,p)] &= \STC H n  &\ & \text{on}\ \Gamma(t), t>0,
  \\\label{eq:4}
V &= n\cdot v|_{\Gamma(t)} - m[n\cdot \nabla\mu] && \text{on}\ \Gamma(t),t>0,\\\label{eq:4'}
\mu|_{\Gamma(t)} & = \STC H&& \text{on}\ \Gamma(t),t>0,\\\label{eq:5}
v|_{\partial\Omega} &= 0 &\ & \text{on}\ \partial\Omega,t>0, \\\label{eq:5'}
 n_\Omega\cdot m\nabla  \mu|_{\partial\Omega}  &= 0 &\ & \text{on}\ \partial\Omega,t>0,\\
\Omega^+(0)&= \Omega_0^+&& \label{eq:6'}\\\label{eq:6}
v|_{t=0} &= v_0 &\ & \text{in}\ \Omega,
\end{alignat}
where $v_0, \Omega_0^+$ are given initial data satisfying
$\partial\Omega_0^+\cap \partial\Omega=\emptyset$ and where $\STC,m>0$ are
a surface tension and a mobility constant, respectively. Implicitly it is
assumed that $v,\mu$ do not jump across $\Gamma(t)$,
i.e., $$[v]=[\mu]=0\qquad \text{on}\ \Gamma(t), t>0.$$  
Equations (\ref{eq:1})-(\ref{eq:2}) describe the conservation of linear
momentum and mass in both fluids and (\ref{eq:3}) is the balance of forces
at the boundary. 
The equations for $v$ are
complemented by the non-slip condition \eqref{eq:5} at the boundary of
$\Omega$. 
The conditions \eqref{eq:2'}, \eqref{eq:5'} describe together with
\eqref{eq:4} a continuity equation for the mass of the phases, and
\eqref{eq:4'} relates the chemical potential $\mu$ to the $L^2$-gradient
of the surface area, which is given by the mean curvature of the
interface. In this formulation of the model we assume (if $\Omega\neq
\mathbb{T}^d$) that $\Gamma(t)$ does not touch $\partial\Omega$.

For $m=0$ the velocity field $v$ is independent of $\mu$. In this case,
(\ref{eq:4}) describes the usual kinematic condition that the 
interface is transported by the flow of the surrounding fluids and
(\ref{eq:1})-(\ref{eq:6}) reduces to the classical model of a two-phase
Navier--Stokes flow as for example studied by cf. Denisova and
Solonnikov~\cite{DenisovaTwoPhase}, where short time existence of strong
solutions is shown. 
On the other hand, if $m>0$, the equations (\ref{eq:2'}), (\ref{eq:4'}), (\ref{eq:5'})
with $v=0$ define the Mullins--Sekerka flow of a family of interfaces.
This evolution describes the gradient flow for the surface area functional
with respect to the $H^{-1}(\Omega)$ scalar product.
Therefore we will also call (\ref{eq:1})-(\ref{eq:6})
the Navier-Stokes/Mullins-Sekerka system. 

The motivation to consider (\ref{eq:1})-(\ref{eq:6}) with $m>0$ is twofold: 
First of all, the modified system gives a regularization of the
classical model $m=0$ since we change from a
parabolic-hyperbolic system to a purely parabolic system (cf. also the
effect of $m>0$ in \eqref{eq:NSCH3} below). 
Secondly, (\ref{eq:1})-(\ref{eq:6}) appears as sharp
interface limit of the following diffuse interface model, introduced by
Hohenberg and Halperin~\cite{HohenbergHalperin} and rigorously derived by
Gurtin et al.~\cite{GurtinTwoPhase}:
\begin{alignat}{2}\label{eq:NSCH1}
  \partial_t v + v\cdot \nabla v - \dive (\nu (c)Dv) + \nabla p &= -\eps\dive
  (\nabla c \otimes \nabla c) 
  &\qquad&\text{in}\ \Omega\times (0,\infty), \\\label{eq:NSCH2}
  \dive v &=0 &\qquad&\text{in}\ \Omega\times (0,\infty), \\\label{eq:NSCH3}
  \partial_t c + v\cdot\nabla c &= m\Delta \mu &\qquad&\text{in}\ \Omega\times (0,\infty), \\\label{eq:NSCH4}
\mu &= \eps^{-1}f'(c) - \eps\Delta c&\qquad& \text{in}\ \Omega\times (0,\infty),\\\label{eq:NSCH5}
  v|_{\partial\Omega} &=0 &\qquad& \text{on}\ \partial\Omega\times (0,\infty),\\\label{eq:NSCH6}
  \partial_n c|_{\partial\Omega} = \partial_n \mu|_{\partial\Omega} &= 0
  &\qquad& \text{on}\ \partial\Omega\times (0,\infty),\\\label{eq:NSCH7}
  (v,c)|_{t=0} &= (v_0,c_0) &\qquad& \text{in}\ \Omega.
\end{alignat}
Here $c$ is the concentration of one of the fluids, where we note that a
partial mixing of both fluids is assumed in the model, and $f$ is a suitable
``double-well potential'' e.g. $f(c)=c^2(1-c)^2$. Moreover,
$\eps>0$ is a small parameter related to the interface thickness, $\mu$
is the so-called chemical potential and $m>0$ the mobility. We refer to
\cite{ModelH} for a further discussion of this model and to the appendix
where we discuss the convergence of 
(\ref{eq:NSCH1})-(\ref{eq:NSCH7}) to varifold solutions of 
(\ref{eq:1})-(\ref{eq:6}).


Sufficiently smooth solutions of (\ref{eq:1})-(\ref{eq:6}) satisfy the
following energy equality,
\begin{equation}\label{eq:en_id_nsms}
  \frac{d}{dt} \frac12 \int_\Omega |v(t)|^2\sd x +\STC \frac{d}{dt}
  \Ha^{d-1}(\Gamma(t)) = - \int_\Omega \nu(\Chi) |Dv|^2 \sd x - m\int_\Omega |\nabla \mu|^2 \sd x, 
\end{equation}
where $\nu(0)= \nu_-$ and $\nu(1) = \nu_+$ and $\Ha^{d-1}$ denotes the
$(d-1)$-dimensional Hausdorff measure. This identity can be
verified by multiplying (\ref{eq:1}) and (\ref{eq:2'}) with $v$, $\mu$,
resp., integrating and using the boundary and interface conditions
(\ref{eq:3})-(\ref{eq:5'}). This energy equality motivates the choice of
solution spaces in our weak formulation and shows that the
regularization introduced for $m>0$ yields an additional dissipation
term. In particular, we expect $\cp(\cdot,t)\in H^{1,2}(\Omega)$ for almost all
$t\in\R_+$ and formally, using Sobolev
inequality and \eqref{eq:5}, that $H(\cdot,t)\in L^4(\Gamma(t))$ for
$d\leq 3$. This gives some indication of extra regularity properties of the
phase interfaces in the model with $m>0$.

Our main result is the existence of weak solutions of
(\ref{eq:1})-(\ref{eq:6}) for large times. For the definitions of the
function spaces we refer to Section~\ref{sec:Notation} below; the
concept of generalized mean curvature for non-smooth phase interfaces is
taken from \cite{Roeg04}, see Definition \ref{def:gen-mc} below. 
\begin{theorem}\label{thm:Main}
Let $d=2,3$, $T>0$, let $\Omega\subseteq \R^d$ be a bounded domain with
smooth boundary or let $\Omega= \mathbb{T}^d$ and set
$\Omega_T=\Omega\times (0,T)$. Moreover, let $\nu(0):=\nu_-$, $\nu(1):=\nu_+$
and $\STC,m>0$. Then for any $v_0\in L^2_\sigma(\Omega)$, $\Chi_0\in
BV(\Omega;\{0,1\})$ there are $v\in
L^\infty(0,T;L^2_\sigma(\Omega))\cap
L^2(0,T;H^1_0(\Omega)^d)$, 
$\Chi\in L^\infty(0,T;BV(\Omega;\{0,1\}))$,  $\mu\in
L^2(0,T;H^1(\Omega))$, that satisfy
(\ref{eq:1})-(\ref{eq:6}) in the following sense: For almost all
$t\in (0,T)$ the phase interface $\partial^*\{\Chi(\cdot,t)=1\}$ has
a generalized mean curvature vector $H(t)\in L^s(d|\nabla\Chi^h(t)|)^d$ 
with $s=4$ if $d=3$ and $1\leq s<\infty$ arbitrary if $d=2$, such that
\begin{eqnarray}\nonumber
\lefteqn{\int_{\Omega_T}\left(-v\partial_t \varphi  +  v\cdot\nabla v \varphi
  +\nu(\Chi)Dv:D\varphi\right)  \,d(x,t)}\\ 
&&-
  \int_\Omega \varphi|_{t=0}\cdot v_0\,dx 
   =\, \STC \int_0^\infty\int_{\Omega} H(t)\cdot \varphi(t)\,d|\nabla \Chi(t)|\,dt, \label{eq:ns}
\end{eqnarray}
holds for all $\varphi \in C^\infty([0,T];
  C_{0,\sigma}^\infty(\Omega)) $ with $\varphi|_{t=T}=0$,
\begin{align}
  \int_{\Omega_T}\Chi(\partial_t\psi + \dive(\psi v)) \,dx\,dt +
  \int_\Omega \Chi_0(x)\psi(0,x)\,dx =  &m\int_{\Omega_T} \nabla\cp\cdot\nabla\psi \,dx\,dt\label{eq:mu} 
\end{align}
holds for all $\psi\in
C^\infty([0,T]\times \overline{\Omega})$ with $\psi|_{t=T}=0$ and 
\begin{equation}\label{eq:HEqn}
  \kappa H(t,.) = \mu(t,.)\frac{\nabla\Chi(\cdot,t)}{|\nabla\Chi|(\cdot,t)} \qquad \Ha^{d-1}-\text{almost everywhere on}\ \partial^\ast \{\Chi(t,.)=1\}
\end{equation}
holds for almost all $0<t<T$.
\end{theorem}

\begin{remark} 
(\ref{eq:ns}) is the weak formulation of (\ref{eq:1}),
(\ref{eq:3}), and \eqref{eq:6}. It is obtained from testing (\ref{eq:1})
  with $\varphi$ in $\Omega^\pm(t)$, integrating over
  $\Omega^+(t)\cup\Omega^-(t)$ and using 
(\ref{eq:4}) together with Gauss' theorem. Similarly, \eqref{eq:mu} is a
weak formulation of \eqref{eq:2'}, \eqref{eq:4}, \eqref{eq:5'}, and
\eqref{eq:6'}. The conditions \eqref{eq:2}, \eqref{eq:5} and $[v]=[\mu]=0$ on 
$\Gamma(t)$ are included in the choice of the function spaces, namely
$v(t)\in H^1_0(\Omega),\mu(t)\in H^1(\Omega)$ for almost every $t\geq
0$, and \eqref{eq:4'} is formulated in \eqref{eq:HEqn}.

Finally, we note that, because of
(\ref{eq:HEqn}), (\ref{eq:ns}) is equivalent to
  \begin{eqnarray}\nonumber
 \lefteqn{
    \int_{\Omega_T}\left(-v\partial_t \varphi  +  v\cdot\nabla v \varphi
    +\nu(\Chi)Dv:D\varphi\right)  \,d(x,t) }\\   \label{eq:ns2}
    &&-\int_\Omega v_0\cdot \varphi|_{t=0}\,dx 
    =\, -\int_{\Omega_T} \Chi\nabla \mu \cdot \varphi\,d(x,t) 
  \end{eqnarray}
for all $\varphi \in C^\infty([0,T];
  C_{0,\sigma}^\infty(\Omega))$ with $\varphi|_{t=T}=0$. The latter form will be used for the
  construction of weak solutions. 
\end{remark}

\begin{remark}
Compared to the available long-time existence results for the classical
model $m=0$ and as a consequence of the diffusive effects that are
included in our model, 
Theorem \ref{thm:Main} yields the long-time existence of  more
regular solutions. In the classical case $m=0$
Plotnikov~\cite{PlotnikovTwoPhase} and
Abels~\cite{GeneralTwoPhaseFlow,ReviewGeneralTwoPhaseFlow} have shown 
the long-time existence of generalized solutions. However, in their
formulations the phase interfaces are in general not regular enough to
define a mean curvature.
Condition \eqref{eq:3} is only satisfied by a varifold
that may depend on the construction process of the solutions, 
that is in general not rectifiable, and lacks a $(d-1)$-dimensional
character  (due to concentration and 
oscillation effects of the interface and  e.g. the formation of
``infinitesimal small droplets'', cf. the discussion in 
\cite{ReviewGeneralTwoPhaseFlow}). In contrast, in our
weak formulation the phase interfaces have a generalized mean curvature
that enjoys the integrability property that are expected, in
the smooth case, from \eqref{eq:5}, the energy equality, and the
Sobolev inequality for the chemical potential. A similar result for the
case $m=0$ is an open problem.

We note that a similar but different regularization was proposed by Liu and
Shen~\cite{LiuShenModelH}. 
In their model (\ref{eq:4}) is replaced by
\begin{equation*}
  V= n\cdot v|_{\Gamma(t)} + m H.
\end{equation*}
Local in time well-posedness for the latter system was proved by
Maekawa~\cite{MaekawaTwoPhaseFlow}. Physically, this model has the
disadvantage that the mass of the fluids, i.e., $|\Omega^\pm(t)|$, is not
preserved in time, while this is the case for our system
(\ref{eq:1})-(\ref{eq:6}).
\end{remark}

\begin{remark}
We note that our concept of weak solution does not include a 
formulation of a contact angle
condition in the case that $\Omega$ is a bounded domain and the phase
boundary $\partial^\ast 
\{\Chi(t,.)=1\}$ meets the boundary of the domain $\partial\Omega$. Even for
weak solutions of the Mullins-Sekerka flow as constructed in
\cite{Roeg04} the formulation of boundary
conditions is an open problem. 
\end{remark}

For simplicity we will assume $\STC=m=1$ in the following. All statement
below and the proof of Theorem~\ref{thm:Main} will be valid for general
$m,\kappa>0$ if modified accordingly. The structure of the article is as
follows:
First the basic notation and some preliminaries are summarized in
Section~\ref{sec:Notation}. Then weak solutions of a time-discrete approximate
system are constructed in Section~\ref{sec:td}. Our main theorem is
proved in Section~\ref{sec:limit} by passing to the limit in the
approximate system. Finally, in the appendix we prove the convergence of
the diffuse interface model \eqref{eq:NSCH1}-\eqref{eq:NSCH7} to
\eqref{eq:1}-\eqref{eq:6}. However, in this limit we have to work with a
weaker notion of generalized solutions, compared to the notion of
solutions that we use in Theorem \ref{thm:Main}. 


\section{Notation and Preliminaries}\label{sec:Notation}


For $A,B\in \R^{d\times d}$ we denote $A:B=\operatorname{tr} (AB)$ and
$|A|=\sqrt{A:A}$. Given $a\in\R^d$
  we define $a\otimes a \in\R^{d\times d}$ as the matrix with the entries
  $a_ia_j$, $i,j=1,\ldots, d$. The space of all $k$-dimensional unoriented
  linear subspaces of $\R^d$ is denoted by $G_k$. If $X$ is a Banach space,
  $X^\ast$ denotes its dual and $\langle x^\ast,x\rangle\equiv \langle x^\ast,
  x\rangle_{X^\ast,X}$, $x^\ast\in X^\ast, x\in X$, the duality product. If $H$
  is a Hilbert space, then $(\cdot,\cdot)_H$ denotes its inner product. Moreover, we
  use the abbreviation $(\cdot,\cdot)_M=(\cdot,\cdot)_{L^2(M)}$.

For $s>0$ we denote by $[s]$ the integer part of $s$ and
for $f:\R\to X$ we define the backward and forward difference quotients
by
\begin{gather*}
  \partial^-_{t,h}f\,:=\, \frac{f(t)-f(t-h)}{h},\quad
  \partial^+_{t,h}f\,:=\, \frac{f(t+h)-f(t)}{h}.
\end{gather*}

\subsection{Measures and BV-Functions}\label{subsec:BVVarifolds}
Let $X$ be a locally compact separable metric space and let $C_0(X;\R^m)$ by the closure of compactly supported continuous functions $f\colon X\to \R^m$, $m\in\N$, in the supremum norm. Moreover, denote by $\M(X;\R^m)$ the space of all finite $\R^m$-valued Radon measures, $\M(X):= \M(X;\R)$, and
$\M_1(X)$ denotes the space of all probability measures on $X$.
By the Riesz representation theorem $\M(X;\R^m)=C_0(X;\R^m)^*$, cf. e.g.
\cite[Theorem 1.54]{AmbrosioEtAl}. Given $\lambda\in \M(X;\R^m)$ we denote by
$|\lambda|$ the total variation measure defined by 
\begin{equation*}
  |\lambda|(A) = \sup \left\{\sum_{k=0}^\infty |\lambda(A_k)|: A_k\in \mathcal{B}(X)\ \text{pairwise disjoint}, A= \bigcup_{k=0}^\infty A_k \right\}
\end{equation*}
for every $A\in \mathcal{B}(X)$, where $\mathcal{B}(X)$ denotes the
$\sigma$-algebra of Borel sets of $X$.  Moreover,
$\frac{\lambda}{|\lambda|}\colon X\to \R^m$
denotes the Radon-Nikodym derivative of $\lambda$ with respect to $|\lambda|$.
The restriction of a measure $\mu$
to a $\mu$-measurable set $A$ is denoted by $(\mu\lfloor A)(B)= \mu(A\cap B)$.
Finally, the $s$-dimensional Hausdorff measure on $\R^d$, $0\leq s\leq d$, is
denoted by $\Ha^s$.

Let $U\subseteq \R^d$ be an open set. 
Recall that 
\begin{eqnarray*}
  BV(U)&=& \{f \in L^1(U): \nabla f \in \M(U;\R^d)\}\\
  \|f\|_{BV(U)}&=& \|f\|_{L^1(U)}+ \|\nabla f\|_{\M(U;\R^d)},
\end{eqnarray*}
where $\nabla f$ denotes the distributional derivative. Moreover,
$BV(U;\{0,1\})$ denotes the set of all $\Chi\in BV(U)$ such that
$\Chi(x)\in\{0,1\}$ for almost all $x\in U$.


A set $E\subseteq U$ is said to have finite perimeter in $U$ if
$\Chi_E\in BV(U)$. By the structure theorem of sets of finite perimeter
$|\nabla \Chi_E|=\Ha^{d-1}\lfloor \partial^\ast E$, where
$\partial^\ast E$ is the so-called reduced boundary of $E$ and for all
$\varphi\in C_0(U,\Rn)$ 
\begin{equation*}
 -\weight{\nabla \Chi_E,\varphi} = \int_E \Div \varphi \sd x =
 -\int_{\partial^\ast E} \varphi\cdot n_E \sd \Ha^{d-1}, 
\end{equation*}
where $n_E(x) = \frac{\nabla \Chi_E}{|\nabla \Chi_E|}$,
cf. e.g. \cite{AmbrosioEtAl}. Note that, if $E$ is a domain with
$C^1$-boundary, then $\partial^\ast E = \partial E$ and $n_E$ coincides
with the interior unit normal.  

\subsection{Function Spaces}
As usual the space of smooth and compactly supported functions in an open set
$U$ is denoted by $C_0^\infty(U)$. Moreover, $C^\infty(\ol{U})$ denotes the
set of all smooth functions with continuous derivatives on $\ol{U}$.   If $X$ is a Banach-space, the $X$-valued
variants are denoted by $C_0^\infty(U;X)$ and
$C^\infty(\ol{U};X)$. For $0<T\leq \infty$, we denote by
$L^p_{loc}([0,T);X)$, $1\leq p\leq \infty$,
the space of all strongly measurable $f\colon (0,T)\to X $ such that $f\in
L^p(0,T';X)$ for all $0<T'<T$.

Furthermore, 
$C_{0,\sigma}^\infty(\Omega)= \{\varphi\in C_0^\infty(\Omega)^d:
\Div \varphi =0\}$ and
\begin{equation*}
  L^2_\sigma(\Omega)\,:=\, \overline{C_{0,\sigma}^\infty(\Omega)}^{L^2(\Omega)},
\end{equation*}
cf. e.g. \cite{Sohr01} for other characterizations and properties of
$L^2_\sigma(\Omega)$.

Finally, we will use the following notation:
\begin{align*}
 \Xmu \,&:=\, H^1(\Omega)\cap \{\cp: \int_\Omega\cp\,=\,
 0\},\\
  \Xns \,&:=\, H^1_0(\Omega,\Rn)\cap L^2_\sigma(\Omega),\\
 \Xms \,&:=\, BV(\Omega;\{0,1\})\cap \left\{\int_\Omega \Chi = m_0\right\},\\
 H^{-1}_{(0)}(\Omega)\,&:=\, H^{1}_{(0)}(\Omega)^*.
\end{align*}
Here we note that $H^1_{(0)}(\Omega)$ is equipped with
the norm $\|f\|_{H^1_{(0)}(\Omega)}=\|\nabla f\|_{L^2(\Omega)}$ and
$H^{-1}_{(0)}(\Omega)$ with the dual norm associated to the latter norm. In
particular, this yields the useful relation
\begin{equation}
  \label{eq:HNegNorm}
  \|f\|_{H^{-1}_{(0)}(\Omega)}= \|\nabla (-\Delta_N)^{-1} f\|_{L^2(\Omega)}
\qquad \text{for all}\ f\in H^{-1}_{(0)}(\Omega),
\end{equation}
where $-\Delta_N:\, \Xmu\,\to\,\Xmustar$ is the weak Laplace operator with Neumann boundary
conditions defined by 
\begin{gather}
\dprod[\Xmustar,\Xmu]{-\Delta_N w,\varphi} =  \int_\Omega \nabla
w\cdot\nabla\varphi\sd x  \quad\text{
  for all }\varphi\in \Xmu.  \label{eq:def-DeltaN}
\end{gather}

\label{sec:intro}

\section{Time-discrete approximation}
\label{sec:td}
In this section we will construct weak solutions of an approximate
time-discrete system. Fortunately the coupling of the Navier-Stokes to the
Mullins-Sekerka system can be treated explicitely. The main
result of this section is:
\begin{proposition}\label{prop:td}
Let the assumptions of Theorem~\ref{thm:Main} be valid and let $m_0=
\int_\Omega \Chi_0\sd x$.
Then for all $h>0$ sufficiently small there exist time-discrete solutions
$v^h\in L^\infty(-h,T;\Xns)$, $\Chi^h\in L^\infty(-h,T;\Xms)$, $\cp^h_0\in
L^\infty(0,\infty;\Xmu)$, and  
Lagrange multipliers $\lambda^h:[0,T]\to \R$ such that for all
$t\in (0,T)$ the following equations hold 
\begin{eqnarray}\nonumber
  \lefteqn{\int_\Omega \Big(\frac{1}{h}\big(v^h(t)-v^h(t-h)\big)\eta + v^h(t-h)\cdot
  \nabla 
  v^h(t) \eta\Big)\sd x}\\
&& +\int_\Omega \nu(\Chi^h(t))Dv^h(t): D\eta \,dx 
  =\, -\int_\Omega \Chi^h(t)\nabla \cp^h_0(t)\cdot \eta\,dx
  \label{eq:td-ns}
\end{eqnarray}
for all $\eta\in \Xns$,
\begin{gather}
  v^h(t)\,=\,v_0,\quad \Chi^h(t)\,=\, \Chi_0 \label{eq:td-init}
\end{gather}
for $-h\leq  t\leq 0$,
\begin{gather}\label{eq:td-mu}
  \int_\Omega \left(\frac{\Chi^h(t)-\Chi^h(t-h)}h \xi -
   (v^h\Chi^h)(t-h)\cdot  \nabla\xi\right)\,dx\,=\, -\int_\Omega
 \nabla\cp^h_0(t)\cdot \nabla \xi \sd x
\end{gather}
for all $\xi\in H^1(\Omega)$, and with $\cp^h(t):= \cp^h_0(t)+\lambda^h(t)$
\begin{gather}\label{eq:td-gt}
  \int_\Omega \Big(\dive\eta - \frac{\nabla \Chi^h(t)}{|\nabla
  \Chi^h(t)|}\cdot D\eta \frac{\nabla \Chi^h(t)}{|\nabla
  \Chi^h(t)|}\Big) \,d|\nabla \Chi^h(t)| \,=\, \int_\Omega \Chi^h(t)\dive
  \big(\cp^h(t)\eta\big)
\end{gather}
for all $\eta\in C^1(\overline{\Omega};\Rn)$ with $\eta\cdot n_\Omega=0$ on $\partial\Omega$.
Moreover, we have the estimates
\begin{eqnarray}\nonumber
  \lefteqn{ \sup_{t\in (0,T)}\frac{1}{2}\|v^h(t)\|_{L^2(\Omega)}^2 +
    \sup_{t\in (0,T)}\int_\Omega \,d|\nabla\Chi^h(t)|}\\
&& +
  \nu_{\min}\|Dv^h\|_{L^2(\Omega_T)}^2
  + \frac{1}{4}\|\nabla\cp^h_0\|_{L^2(\Omega_T)}^2 
  \leq\, C\Big(\|v_0\|_{L^2(\Omega)},\int_\Omega \,d|\nabla\Chi_0|\Big)
  \label{eq:est-td}
\end{eqnarray}
and for all $t\in (0,T)$
\begin{gather}
  \|\cp^h_0(t)+\lambda^h(t)\|_{H^1(\Omega)} \,\leq\,
  C\Big(1+\|\nabla\cp^h_0(t)\|_{L^2(\Omega)}\Big) \label{eq:est-lagrange}
\end{gather}
holds, where $C$ depends only on $d$, $\Omega$, $T$, $m_0$, the initial data,
and $\nu_{min}=\min (\nu(0),\nu(1))$.
\end{proposition}
In the remainder of this section we prove Proposition \ref{prop:td}. The
first step gives the solvability and estimates for a time-discrete
and regularized Navier--Stokes equation.
\begin{lemma}\label{lem:td-ns}
Let $\tilde{v}\in \Xns$, $\Chi\in \Xms$, $\cp\in H^1(\Omega)$, and $h>0$ be given. 
Then there exists a solution $v\in \Xns$ of
\begin{gather}\label{eq:lem-td-ns}
  \int_\Omega \left(\frac{1}{h}\left(v-\tilde{v}\right) +\tilde{v}\cdot\nabla
    v\right) \varphi \, dx
  + \int_\Omega \nu(\Chi)Dv:D\varphi\, dx \,=\,
  -\int_\Omega \Chi\nabla \cp\cdot \varphi \, dx
\end{gather}
for all $\varphi\in C^\infty_{0,\sigma}(\Omega)$.
For each solution $v\in \Xns$ 
\begin{gather}\label{eq:lem-est-td-ns}
  \frac{1}{2}\|v\|_{L^2(\Omega)}^2 + h\int_\Omega \nu(\Chi)|Dv|^2
  \,\leq\, \frac{1}{2}\|\tilde{v}\|_{L^2(\Omega)}^2 - h
  \int_\Omega \Chi\nabla\cp\cdot v\,dx
\end{gather}
holds.
\end{lemma}
\begin{proof}
We show the existence of a solution $v\in \Xns$ of (\ref{eq:lem-td-ns}) that
satisfies (\ref{eq:lem-est-td-ns}) with the aid of the Leray-Schauder
principle, cf. e.g. \cite[Chapter~II, Lemma~3.1.1]{Sohr01}. To this end, we
define $L\colon \Xns \to \Xns^\ast$ and $G\colon \Xns\to L^3_\sigma(\Omega)^\ast\cong L^{3/2}_\sigma(\Omega)$ by
\begin{eqnarray*}
  \left\langle Lv,\varphi\right\rangle_{\Xns',\Xns} &=& \int_\Omega
  \nu(\Chi)Dv:D\varphi \sd x
  \\
  \left( G(v),\psi\right)_\Omega &=& \int_\Omega \left(-\Chi\nabla\cp
    -\frac1h(v-\tilde{v}) -\tilde{v}\cdot \nabla v\right)\cdot \psi \, dx 
\end{eqnarray*}
for all $v,\varphi \in \Xns$, $\psi \in L^3_\sigma(\Omega)$.
By the Lemma of Lax-Milgram, $L\colon \Xns\to
\Xns'$ is an isomorphism. Moreover, it is easy to check that $G\colon \Xns\to
L^3_\sigma(\Omega)^\ast$ is a continuous mapping, where we note that
$\tilde{v}\cdot \nabla v\in L^{3/2}(\Omega)$. Since $L^3_\sigma(\Omega)^\ast\hookrightarrow
\Xns^\ast$ compactly, $G\colon \Xns\to \Xns^\ast$ is completely continuous. Thus $F
= L^{-1} G\colon \Xns\to \Xns$ is completely continuous and
(\ref{eq:lem-td-ns}) is equivalent to the fixed-point problem
$
  v= G(v)
$.
In order to apply the Leray-Schauder principle, let  $R>0$ be such that $R^2= M^2\|\nabla
  \cp\|_{L^2(\Omega)}^2+ \frac{M}h\|\tilde{v}\|_{L^2(\Omega)}^2$, where $M$
is a constant such that 
\begin{equation*}
  \|v\|_{H^1(\Omega)}^2 \leq M\int_\Omega \nu(\Chi) |Dv|^2 dx\qquad \text{for
    all}\ v\in \Xns.
\end{equation*}
It remains to show that for all $v\in \Xns$
\begin{equation}
  \label{eq:LSCond}
  v= \lambda F(v), \lambda\in[0,1] \quad \Rightarrow \quad \|v\|_{\Xns}\leq R. 
\end{equation}
To this end let $v= \lambda F(v)$ for some $\lambda\in [0,1]$. Then $Lv=
\lambda G(v)$ and therefore
\begin{equation*}
    \frac{\lambda}{h}\int_\Omega \left(v-\tilde{v}\right)\varphi\, dx  +
    \int_\Omega \lambda \tilde{v}\cdot\nabla v\varphi\, dx 
  +\int_\Omega \nu(\Chi)Dv:D\varphi\, dx \,=\,
  \lambda\int_\Omega \Chi\nabla \cp \cdot\varphi\, dx
\end{equation*}
for all $\varphi \in \Xns$.
Choosing $\varphi = v$ yields
\begin{equation}\label{eq:vEstim}
  \frac{\lambda}{2h}\|v\|_{L^2(\Omega)}^2 +\int_\Omega \nu(\Chi)|Dv|^2\,
  dx \leq 
  \lambda\int_\Omega \Chi\nabla\cp\cdot v\,dx +
  \frac{\lambda}{2h}\|\tilde{v}\|_{L^2(\Omega)}^2, 
\end{equation}
where we have used $(v-\tilde{v})\cdot v= \frac12|v|^2-\frac12
|\tilde{v}|^2+\frac12 |v-\tilde{v}|^2\geq \frac12|v|^2-\frac12
|\tilde{v}|^2$. Hence, using that $|\int_\Omega \Chi\nabla\cp\cdot
v|\leq \|\nabla \cp\|_{L^2(\Omega)}\|v\|_{L^2(\Omega)}$, 
\begin{eqnarray*}
  \|v\|_{H^1(\Omega)}^2&\leq& M \int_\Omega \nu(\Chi)Dv:Dv\,
  dx \leq M\left(\|\nabla \cp\|_{L^2(\Omega)}\|v\|_{L^2(\Omega)} +
    \frac1{2h}\|\tilde{v}\|_{L^2(\Omega)}^2\right)\\
  &\leq &\frac{M^2}2 \|\nabla \cp\|_{L^2(\Omega)}^2 +
    \frac{M}{2h}\|\tilde{v}\|_{L^2(\Omega)}^2+ \frac12\|v\|_{H^1(\Omega)}^2
\end{eqnarray*}
and therefore $\|v\|_{H^1(\Omega)}\leq R$.

Because of (\ref{eq:LSCond}), the Leray-Schauder principle implies the
existence of a fixed point $v\in \Xns$ which solves
(\ref{eq:lem-td-ns}). Finally, (\ref{eq:lem-est-td-ns}) follows from
(\ref{eq:vEstim}) with $\lambda =1$. 
\end{proof}
Next we solve the appropriate versions of the Mullins--Sekerka part
\eqref{eq:td-mu}, \eqref{eq:td-gt}. We follow
Luckhaus--Sturzenhecker \cite{LStu95} and use that the Mullins--Sekerka
flow is the $H^{-1}$-gradient flow of the surface-area-functional.
\begin{lemma}\label{lem:td-ms}
For $\tilde{\Chi}\in \Xms$ and $\tilde{v}\in \Xns $ there
exist $\Chi\in \Xms$, $\cp_0\in \Xmu$, and a constant $\lambda\in\R$
such that
\begin{gather}\label{eq:lem-td-mu}
  \int_\Omega \Big(\nabla\cp_0\cdot\nabla\xi + \frac{1}{h}\big(\Chi-\tilde{\Chi}\big)\xi
  -\tilde{v}\tilde{\Chi}\cdot \nabla\xi \Big)\,dx\,=\, 0
\end{gather}
for all $\xi\in H^1(\Omega)$, such that
\begin{gather}\label{eq:lem-td-gt}
  \int_\Omega \Big(\dive\eta - \frac{\nabla \Chi}{|\nabla
  \Chi|}\cdot D\eta \frac{\nabla \Chi}{|\nabla
  \Chi|} \Big)\,d |\nabla \Chi| \,=\, \int_\Omega \Chi\dive
  \big((\cp_0+\lambda)\eta\big)\sd x
\end{gather}
for all $\eta\in C^1(\Omega,\Rn)$ with $\eta\cdot n_\Omega=0$ on
$\partial\Omega$, and such that 
\begin{gather}\label{eq:lem-est-ms}
  \per{\Chi} +
  \frac{h}{2}\|\nabla\cp_0\|_{L^2(\Omega)}^2\,\leq\, \int_\Omega
  d|\nabla\tilde{\Chi}| +\frac{h}{2}\|\tilde{v}\|_{L^2(\Omega)}^2.
\end{gather}
Moreover, we have
\begin{gather}
  |\lambda|\,\leq\, C(m_0,n,\Omega)\Big(1+\per{\Chi}\Big)
  \Big(\per{\Chi} + \|\nabla\cp_0\|_{L^2(\Omega)}\Big),
  \label{eq:est-lambda-1}\\
  \|\cp_0+\lambda\|_{H^1(\Omega)}\,\leq\, C(m_0,n,\Omega)\Big(1+\per{\Chi}\Big)
  \Big(\per{\Chi} + \|\nabla\cp_0\|_{L^2(\Omega)}\Big).
  \label{eq:est-lambda} 
\end{gather}
\end{lemma}
\begin{proof}
We split the proof into several steps.
\step{1}
There exist $\cp_0\in\Xmu, \Chi\in\Xms$ satisfying \eqref{eq:lem-td-mu},
\eqref{eq:lem-est-ms}, and enjoying a minimizing property
from which we will deduce \eqref{eq:lem-td-gt}.\\
Define a functional $F^h:\Xms \to\R$,
\begin{gather}\label{eq:def-Fh} 
  F^h(\sigma)\,=\, \int_{\Omega}|\nabla \sigma| +
  \frac{1}{2h}
  \|\sigma-\tilde{\Chi}+h\tilde{v}\cdot\nabla\tilde{\Chi}\|_{H^{-1}_{(0)}(\Omega)}^2
\end{gather}
for $\sigma\in \Xms$. 
We remark that $\tilde{v}\cdot\nabla\tilde{\Chi}\in\Xmustar$ is defined by
\begin{gather}
  \dprod{\tilde{v}\cdot\nabla\tilde{\Chi},\zeta}\,=\, -\int_\Omega
  \tilde{\Chi}\tilde{v}\cdot\nabla\zeta\sd x, \label{eq:def-dprod-vChi}
\end{gather}
where we note that $\Div \tilde{v}=0$.
Because of (\ref{eq:def-DeltaN}),
\begin{eqnarray}\nonumber
  \lefteqn{\|\sigma-\tilde{\Chi}+h\tilde{v}\cdot\nabla\tilde{\Chi}\|_{H^{-1}_{(0)}(\Omega)}^2}\\
  &=&
  \bdprod[H^1_{(0)}(\Omega)^\ast,H^1_{(0)}(\Omega)]{\sigma-\tilde{\Chi}+h\tilde{v}\cdot\nabla\tilde{\Chi},  
    (-\Delta_N)^{-1}\big(\sigma-\tilde{\Chi}+h\tilde{v}\cdot\nabla\tilde{\Chi}\big)}.
\end{eqnarray}
By the $L^1(\Omega)$-compactness of bounded sequences in $BV(\Omega)$,
the lower semi-continuity of the perimeter under
$L^1(\Omega)$-convergence, and the continuity of the embedding
$L^2_{(0)}(\Omega)\to H^{-1}_{(0)}(\Omega)$ there exists a minimizer
$\Chi\in\Xms$ of $F^h$. Moreover, 
\begin{gather}\label{eq:def-td-cp}
  \cp_0\,:=\, -(-\Delta_N)^{-1}\Big(\frac{1}{h}\big(\Chi-\tilde{\Chi}\big)
  + \tilde{v}\cdot\nabla\tilde{\Chi}\Big) 
\end{gather}
satisfies \eqref{eq:lem-td-mu}. 
We deduce now \eqref{eq:lem-est-ms} from
$F^h(\Chi)-F^h(\tilde{\Chi})\leq 0$ and \eqref{eq:def-DeltaN},
\eqref{eq:def-td-cp}. In fact,
\begin{align*}
  \per{\Chi} +\frac{1}{2h}\int_\Omega |h\nabla\cp_0|^2\, dx \,&\leq\,
  \per{\tilde{\Chi}} +\frac{h}{2}
  \left\langle\tilde{v}\cdot\nabla\tilde{\Chi}, 
    (-\Delta_N)^{-1}\big(\tilde{v}\cdot\nabla\tilde{\Chi}\big)\right\rangle_{H^{-1}_{(0)},H^1_{(0)}}  \\ 
  &=\,\per{\tilde{\Chi}} +\frac{h}{2}\int_\Omega
  |\nabla  (-\Delta_N)^{-1}\tilde{v}\cdot\nabla\tilde{\Chi}|^2 \sd x\\
  &\leq\,  \per{\tilde{\Chi}}
  +\frac{h}{2}\|\tilde{v}\|_{L^2(\Omega)}^2,
\end{align*}
where in the last step we have used that by \eqref{eq:def-DeltaN},
\eqref{eq:def-dprod-vChi} 
\begin{align*}
  \|\nabla
  (-\Delta_N)^{-1}\tilde{v}\cdot\nabla\tilde{\Chi}\|_{L^2(\Omega)}^2 \,&=\, 
  \left\langle
  \tilde{v}\cdot\nabla\tilde{\Chi},
  (-\Delta_N)^{-1}\tilde{v}\cdot\nabla\tilde{\Chi}\right\rangle_{H^{-1}_{(0)},H^1_{(0)}}\\
  &=\, -\int_\Omega \tilde{\Chi}\tilde{v}\cdot \nabla
  (-\Delta_N)^{-1}\tilde{v}\cdot\nabla\tilde{\Chi}\sd x\\
  &\leq\, \|\tilde{v}\|_{L^2(\Omega)}\|\nabla
  (-\Delta_N)^{-1}\tilde{v}\cdot\nabla\tilde{\Chi}\|_{L^2(\Omega)},
\end{align*}
hence $\|\nabla  (-\Delta_N)^{-1}\tilde{v}\cdot\nabla\tilde{\Chi}\|_{L^2(\Omega)}
  \,\leq\,\|\tilde{v}\|_{L^2(\Omega)}$.
\step{2} We compute the first variation of $F^h$ in $\Chi$ with respect to volume
preserving variations. With this 
aim we consider a smooth family $(\Phi_s)_{s\in (-\eps,\eps)}$ of smooth
diffeomorphisms $\Phi_s:\Omega\to\Omega$ with $\Phi_0=\Id$ and variation
field $\eta$ such that
\begin{gather*}
  \eta\,=\, \left.\frac{\partial}{\partial s}\right|_{s=0}\Phi_s\in C^\infty(\Omega,\Rn),\qquad
  \eta\cdot n_\Omega=0\quad \text{ on }\partial\Omega.
\end{gather*}
Assume that the variations $\Phi_s$ conserve the
volume of $\{\Chi=1\}$, that means that $\sigma_s:=
\Chi\circ\Phi_s^{-1}$ satisfy
\begin{gather*}
  \int_\Omega \sigma_s\,dx\,=\, m_0\qquad \text{ for all }-\eps<s<\eps,
\end{gather*}
in particular $\sigma_s\in \Xms$ and $\int_\Omega \Chi\dive\eta =0$.
Since $\Chi$ minimizes $F^h$ in $Z$ we have $\frac{d}{ds}|_{s=0}F^h(\Chi_s)=0$.
The first part of $F^h$ is given by the perimeter-functional
$\Per_\Omega$ and 
\begin{gather}\label{eq:lem-td-ms-1}
 \frac{d}{ds}|_{s=0}\Per_\Omega(\sigma_s)\,=\,
  \delta\Per_\Omega(\Chi)(\eta)\,=\, \int_\Omega \Big(\dive\eta -   
  \frac{\nabla \Chi}{|\nabla 
  \Chi|}\cdot D\eta \frac{\nabla \Chi}{|\nabla
  \Chi|}\Big) \,d|\nabla \Chi|. 
\end{gather}
Let $K^h$ denote the second part of $F^h$.
Since $(-\Delta_N)$ is linear and symmetric, since
$\frac{\partial}{\partial s}|_{s=0}\sigma_s = -\nabla \Chi\cdot\eta$,
and by \eqref{eq:def-td-cp} we obtain that
\begin{gather}\label{eq:lem-td-ms-2}
  \left.\frac{d}{ds}\right|_{s=0}K^h(\sigma_s)\,=\, \delta K^h(\Chi)(\eta)\,=\,
  \left\langle \nabla\Chi\cdot \eta, \cp_0\right\rangle_{H^{-1}_{(0)},H^1_{(0)}}=-\int_\Omega \Chi\dive (\eta \cp_0)\sd x. 
\end{gather}
We therefore deduce from the minimality of $\Chi$ that
\begin{gather}
  0\,=\, \int_\Omega \Big(\dive\eta -
  \frac{\nabla \Chi}{|\nabla 
  \Chi|}\cdot D\eta \frac{\nabla \Chi}{|\nabla
  \Chi|}\Big) \,d|\nabla \Chi| - \int_\Omega \Chi\dive( \eta
  \cp_0)\,dx.   \label{eq:lem-td-gt-1}
\end{gather}
\step{3}
We next prove \eqref{eq:lem-td-gt}. Fix $\xi\in C^\infty_0(\Omega,\Rn)$ such that 
\begin{gather*}
  \int_\Omega \Chi \dive \xi\, dx \,\neq\, 0,\qquad
  \xi\cdot n_\Omega\,=\, 0\quad\text{ on }\partial\Omega,
\end{gather*}
and choose a family $(h_r)_{r\in (-\eps_1,\eps_1)}$ of smooth
diffeomorphisms of $\Omega$ with $h_0\,=\, \Id$ and
$\frac{\partial}{\partial r}|_{r=0}h_r= \xi$. 
Similarly, for given $\eta\in C^\infty(\Omega,\Rn)$ with $\eta\cdot n_\Omega=0$ on
$\partial\Omega$ we let $(g_s)_{s\in (-\eps_1,\eps_1)}$ be a family of smooth
diffeomorphisms of $\Omega$ with $g_0=\Id$ and $\frac{\partial}{\partial
  s}|_{s=0}g_s=\eta$. 
Then the function $f: (-\eps_1,\eps_1)^2\,\to\,\R$ defined by
\begin{gather*}
  f(s,r)\,:=\, \int_\Omega \Chi\circ(g_s\circ h_r)^{-1}\,dx - m_0\,=\,
  \int_\Omega  \left(\det\big((Dg_s\circ h_r)Dh_r\big)-1\right)\,dx 
\end{gather*}
for sufficiently small $\eps_1>0$
satisfies
\begin{gather*}
  f(0,0)\,=\, 0,\qquad \partial_r|_{s=r=0} f(s,r)\,=\, \int_\Omega
  \Chi\dive\xi\,\neq\, 0.
\end{gather*}
Since $f$ is smooth we obtain by the implicit function theorem that
there exists $0<\eps\leq\eps_1$ and
a smooth function $\varrho:(-\eps,\eps)\,\to\, (-\eps_1,\eps_1)$ such
that
\begin{gather}
  f(s,\varrho(s))\,=\, 0 \quad\text{ for all }s\in
  (-\eps,\eps). \label{eq:f-vol}
\end{gather}
We therefore deduce that
\begin{gather*}
  0\,=\, \frac{d}{ds}|_{s=0}f(s,\varrho(s))\,=\,\int_\Omega
  \Big(\dive\eta + \varrho'(0)\dive\xi\Big)\Chi\,dx,
\end{gather*}
hence
\begin{gather}
  \varrho'(0)\,=\, -\Big(\int_\Omega \Chi\dive\xi\,dx\Big)^{-1}
  \int_\Omega \Chi\dive\eta\,dx.
  \label{eq:r-prime}
\end{gather}
By \eqref{eq:f-vol} the family $(g_s\circ h_{r(s)})_{s\in (-\eps,\eps)}$
defines a variation of $\Omega$ that conserves the volume of $\Chi$.
The corresponding variation field is given by
\begin{gather*}
  \tilde{\eta}\,:=\, \frac{\partial}{\partial s}|_{s=0}\big(g_s\circ
  h_r\big)\,=\, \eta +   \varrho'(0)\xi. 
\end{gather*}
Step 2 therefore implies that
\begin{gather*}
  0\,=\, \int_\Omega \Big(\dive\tilde{\eta} -
  \frac{\nabla \Chi}{|\nabla 
  \Chi|}\cdot D\tilde{\eta} \frac{\nabla \Chi}{|\nabla
  \Chi|}\Big) \,d|\nabla \Chi| - \int_\Omega \Chi\dive( \tilde{\eta}
  \cp_0)\,dx  
\end{gather*}
and yields by \eqref{eq:r-prime}
\begin{align}
   &\int_\Omega \Big(\dive\eta -
  \frac{\nabla \Chi}{|\nabla 
  \Chi|}\cdot D\eta \frac{\nabla \Chi}{|\nabla
  \Chi|}\Big) \,d|\nabla \Chi| - \int_\Omega \Chi\dive(\eta
  \cp_0)\,dx  \notag
  \\ &=\, -\varrho'(0)\Big[\int_\Omega \Big(\dive\xi -
  \frac{\nabla \Chi}{|\nabla 
  \Chi|}\cdot D\xi \frac{\nabla \Chi}{|\nabla
  \Chi|}\Big) \,d|\nabla \Chi| - \int_\Omega \Chi\dive(\xi
  \cp_0)\,dx\Big] \notag\\
  &=\, \lambda \int_\Omega\Chi\dive\eta
\end{align}
with
\begin{eqnarray}\nonumber
  \lambda &:=& \Big(\int_\Omega
  \Chi\dive\xi\,dx\Big)^{-1}\\
&&\quad \cdot\Big[\int_\Omega \Big(\dive\xi - 
  \frac{\nabla \Chi}{|\nabla 
  \Chi|}\cdot D\xi \frac{\nabla \Chi}{|\nabla
  \Chi|}\Big) \,d|\nabla \Chi|- \int_\Omega \Chi\dive(\xi
  \cp_0)\,dx\Big]. \label{eq:def-lambda}
\end{eqnarray}
This proves \eqref{eq:lem-td-gt}.
\step{4}
Finally we derive \eqref{eq:est-lambda} by choosing
a particular $\xi$ in Step 3. We adapt the proof given in
\cite{Chen96}. First we choose a Dirac
sequence $(\varphi_\delta)_{\delta>0}$ with kernel $\varphi\in
C^\infty_c(B_1(0))$, $0\leq\varphi\leq 1$, and set
\begin{gather*}
  \Chi_\delta\,:=\, \Chi *\varphi_\delta,\qquad \bar{\Chi}_\delta \,:=\,
  \frac{1}{|\Omega|}\int_\Omega \Chi_\delta.
\end{gather*}
Let $\psi:\Omega\to\R$ be the solution of
\begin{gather*}
  \Delta\psi\,=\, \Chi_\delta -\bar{\Chi}_\delta \quad\text{ in
  }\Omega,\\
  \nabla\psi\cdot n_\Omega\,=\, 0 \quad\text{ on }\partial\Omega,\qquad
  \int_\Omega \psi\,=\,0
\end{gather*}
and choose $\xi:=\nabla\psi$ in Step 3.
We observe that $|\Chi_\delta-\bar{\Chi}_\delta|\,\leq\, 1$ and
$|\nabla \Chi_\delta|\,\leq\,C(\Omega)\delta^{-1}$. By standard
elliptic estimates we conclude that
\begin{gather}
  \|\psi\|_{C^2(\bar{\Omega})}\,\leq\,
  \frac{1}{\delta}C(\Omega). \label{eq:est-psi} 
\end{gather}
Moreover, we compute that
\begin{gather}
  \|\Chi-\Chi_\delta\|_{L^1(\Omega)}\,\leq\,
  C(\Omega)\delta\left(1+\per{\Chi}\right), \label{eq:est-diff-delta}
\end{gather}
and
\begin{gather}
  \bar{\Chi}_\delta \,=\, \frac{1}{|\Omega|}\int_\Omega \Chi_\delta \,\leq\,
  \frac{1}{|\Omega|}\int_{\Rn}\Chi_\delta\,=\, \frac{m_0}{|\Omega|}. \label{eq:est-bar-Chi}
\end{gather}
Therefore we deduce the estimate
\begin{align}
  \int_\Omega \Chi\dive\xi\sd x \,=\,\int_\Omega \Chi\big(\Chi_\delta
  -\bar{\Chi}_\delta\big)\sd x \,&=\, 
  (1-\bar{\Chi}_\delta)m_0 + \int_\Omega \big(\Chi_\delta -\Chi\big)\Chi 
  \notag\\
  &\geq\, \left(1 - \frac{m_0}{|\Omega|}\right)m_0 -
  C(\Omega)\delta\Big(1+\per{\Chi}\Big) \notag\\
  &\geq c(m_0,\Omega), \label{eq:proof-lagrange2}
\end{align}
for $\delta=\delta_0(m_0,\Omega)\big(1+\per{\Chi}\big)^{-1}$.
Further we compute that
\begin{eqnarray}\nonumber
  \lefteqn{\left| \int_\Omega \Big(\dive\xi - 
  \frac{\nabla \Chi}{|\nabla 
  \Chi|}\cdot D\xi \frac{\nabla \Chi}{|\nabla
  \Chi|}\Big) \,d|\nabla \Chi| - \int_\Omega \Chi\dive(\xi
  \cp_0)\,dx\right|}\\\nonumber
  &\leq&
   \|\psi\|_{C^2(\overline{\Omega})}\int_\Omega
   \,d|\nabla\Chi| +
   2\|\cp_0\|_{H^1(\Omega)}\|\psi\|_{C^2(\overline{\Omega})} \\\nonumber
   &\leq &\frac{C(\Omega)}{\delta}\int_\Omega
   \,d|\nabla\Chi| +
   C(n,\Omega)\frac{1}{\delta}\|\nabla\cp_0\|_{L^2(\Omega)} \\
   &\leq&
   C(m_0,n,\Omega)\Big(1+\per{\Chi}\Big)
  \Big(\per{\Chi} + \|\nabla\cp_0\|_{L^2(\Omega)}\Big),
  \label{eq:est-lagrange3} 
\end{eqnarray}
where in the last two steps we have used \eqref{eq:est-psi} and Poincar\'{e}'s
inequality. We now obtain \eqref{eq:est-lambda-1} from
\eqref{eq:def-lambda}, \eqref{eq:proof-lagrange2}, and
\eqref{eq:est-lagrange3}. The estimate \eqref{eq:est-lambda} follows
again from Poincar\'{e}'s inequality.
\end{proof}

\begin{proof}[Proof of Proposition \ref{prop:td}]
We construct iteratively time-discrete solutions $v^h,\Chi^h$. First
set
\begin{gather*}
  v^h(t)\,:=\, v_0, \quad \Chi^h(t)\,:=\, \Chi_0 \qquad\text{ for }
  -h< t\leq 0.
\end{gather*}
Given functions $v^h(t-h), \Chi^h(t-h)$ for $t\in (kh,(k+1)h]\subset [0,T]$, $k\in\N_0$,  Lemma
\ref{lem:td-ms} yields a solution $\Chi^h(t)\in \Xms$, $\cp^h_0(t)\in
\Xmu$, $\lambda^h(t)\in\R$ that satisfy \eqref{eq:td-mu},
\eqref{eq:td-gt} and
\begin{alignat}{1}
&  \per{\Chi^h(t)} +
  \frac{1}{2h}\|\nabla\cp^h_0(t)\|_{L^2(\Omega)}^2\leq
  \per{\Chi^h(t-h)}  +\frac{h}{2}\|v^h(t-h)\|_{L^2(\Omega)}^2,
  \label{eq:proof-est-ms}\\\nonumber 
 &\|\cp^h_0(t)+\lambda^h(t)\|_{H^1(\Omega)}\\ &\ \leq
  C(m_0,n,\Omega)\Big(1+\int_\Omega\,d|\nabla\Chi^h(t)|\Big) 
  \Big(\int_\Omega\,d|\nabla\Chi^h(t)| +
  \|\nabla\cp^h(t)\|_{L^2(\Omega)}\Big). \label{eq:td-gt-rhs}
\end{alignat}
Then we deduce from Lemma \ref{lem:td-ns} the existence of $v^h(t)\in
\Xns$ that satisfies \eqref{eq:td-ns} and the estimate
\begin{align}
  &\frac{1}{2}\|v^h(t)\|_{L^2(\Omega)}^2 + h\int_\Omega
  \nu(\Chi)|Dv^h(t)|^2 \notag\\
  \leq\, &\frac{1}{2}\|v^h(t-h)\|_{L^2(\Omega)}^2 + h
  \|\nabla\cp_0^h(t)\|_{L^2(\Omega)}\|v^h(t)\|_{L^2(\Omega)}\notag\\
  \leq\, &\frac{1}{2}\|v^h(t-h)\|_{L^2(\Omega)}^2 + \frac{h}{4}
  \|\nabla\cp_0^h(t)\|_{L^2(\Omega)}^2 +h\|v^h(t)\|_{L^2(\Omega)}^2.
  \label{eq:proof-est-td-ns} 
\end{align}
By construction $v^h,\Chi^h$ are constant on each
subinterval $(kh,(k+1)h]\subset [0,T]$, $k\in\N$. Summing
\eqref{eq:proof-est-ms}, \eqref{eq:proof-est-td-ns} for $t_k=kh$,
$k=1,...,[t/h]$ we obtain that
\begin{align}
  &\frac{1}{2}\|v^h(t)\|_{L^2(\Omega)}^2 + \int_\Omega \,d|\nabla\Chi^h(t)| +
  \int_0^{h[t/h]}\int_\Omega \Big(\nu_{\min}|Dv^h(\tau)|^2
  + \frac{1}{2}|\nabla\cp_0^h(\tau)|^2\Big)\,dx\,d\tau \notag\\
  \leq\,& \frac{1}{2}\|v_0\|_{L^2(\Omega)}^2 + \int_\Omega
  \,d|\nabla\Chi_0| \notag\\ 
  &+2\int_{-h}^{h[t/h]}\int_\Omega |v^h(\tau)|^2\,dx\,d\tau +
  \frac{1}{4}\int_0^{h[t/h]}\int_\Omega  |\nabla\cp_0^h(t)|^2\,dx\,d\tau.
  \label{eq:td-est-sum}
\end{align}
Using Gronwall's Lemma we deduce \eqref{eq:est-td}. Finally,
\eqref{eq:est-lagrange} follows from \eqref{eq:est-td} and
\eqref{eq:td-gt-rhs}.
\end{proof}

%
\section{Passing time--discrete approximations to a limit}
\label{sec:limit}
We first show strong compactness of $v^h,\Chi^h$. To this end, we will apply
the following theorem by Simon~\cite[Theorem~6]{Simo87}

\begin{theorem}\label{thmsimon} 
Let $X\subset B\subset Y$ be Banach spaces with compact embedding $X\hookrightarrow B$.
Let $T>0$ and let $\mathcal F$ be a bounded subset of $L^q(0,T;X)$, $1<q\leq\infty$.
Assume that for every $0<t_1<t_2<T$
\begin{equation}
  \label{eq:UniformCont}
\sup_{f\in{\mathcal F}}\|\tau_s f -f\|_{L^1(t_1,t_2;Y)} \to 0 \quad \text{as }
s\to 0,
\end{equation}
where $\tau_s f(t):=f(t+s)$ for every $t\in(0,T-s)$.
Then $\mathcal F$ is relatively compact in $L^p(0,T;B)$ for every $p\in[1,q)$.
\end{theorem}
First of all, because of (\ref{eq:td-ns}), $v^h(.-h)\cdot \nabla v^h= \Div
(v^h(.-h)\otimes v^h)$, and since $v^h(.-h)\otimes v^h \in
L^{\frac43}(0,T;L^2(\Omega))$ is bounded, 
$\partial_{t,h}^{-} v^h\in L^{\frac43}(0,T;H^{-1}(\Omega))$, $h\in (0,1)$ is
bounded.
This implies that
\begin{equation}
  \|\tau_{kh} v^h-v^h\|_{L^{\frac43}(0,T-kh;H^{-1})}\leq kh \|\partial_{t,h}^{-}
  v^h\|_{L^{\frac43}(0,T;H^{-1})}\leq Ckh 
\end{equation}
for all $k\in\N$ such that $kh <T$. Therefore 
(\ref{eq:UniformCont}) holds for $\mathcal F=\{v^h:h\in (0,1)\}$ and
$Y=H^{-1}$
since $v^h$ is piecewise constant. Moreover, since $v^h\in
L^4(0,T;H^{\frac12}(\Omega))$ is bounded, 
Theorem~\ref{thmsimon} implies that $(v^h)_{0<h<1}$ is relatively compact in
$L^2(\Omega_T)$. 

Similarly, (\ref{eq:td-mu}) implies that $\partial_{t,h}^- \Chi^h\in
L^2(0,T;H^{-1}(\Omega))$, $h\in (0,1)$ is bounded. Moreover, since $\Chi^h\in
L^\infty(0,T;BV(\Omega))$ is bounded and $BV(\Omega)\hookrightarrow
L^1(\Omega)$ compactly, we obtain that $\Chi^h\in L^1(\Omega_T)$, $0<h<1$, is relatively
compact.
Since $\|\Chi^h\|_{L^\infty(\Omega_T)}=1$, $(\Chi^h)_{0<h<1}$, is relatively
compact in  $L^p(\Omega_T)$ for every $1\leq p<\infty$.

As a corollary we obtain the compactness of time-discrete approximations.
\begin{proposition}\label{prop:comp-td}
Choose a sequence $h\to 0$ and let
$(v^{h},\Chi^{h},\cp^{h},\lambda^{h})$ denote the time-discrete
solutions constructed in Proposition \ref{prop:td}. Then there exists a
subsequence $h_k\to 0$ $(k\to \infty)$ and  $v\in
L^2(0,T;H^1(\Omega)^d)\cap L^\infty(0,T;L^2_\sigma(\Omega))$, 
$\cp\in L^2(0,T;H^1(\Omega))$, $\Chi\in L^\infty(0,T;\Xms)$, and
$\lambda\in L^2(0,T)$
such that
\begin{alignat}{2}
  v_k\,&\to_{k\to\infty}\, v&\qquad&\text{ in } L^2(\Omega_T), \label{eq:conv-v-strong}\\
  v_k\,&\schwto_{k\to\infty}\, v& &\text{ in } L^2(0,T;H^1(\Omega)),
  \label{eq:conv-v-weak}\\ 
  \Chi_k\,&\to_{k\to\infty}\, \Chi &&\text{ in } L^p(\Omega_T)\quad\qquad\text{ for all
  }1\leq p<\infty, \label{eq:conv-Chi} \\ 
  \cp_k \,&\schwto_{k\to\infty}\, \cp &&\text{ in }
  L^2(0,T;H^1(\Omega)). \label{eq:conv-mu}\\
\lambda_k ,&\schwto_{k\to\infty}\, \lambda &&\text{ in }
  L^2(0,T). \label{eq:conv-lambda}
\end{alignat}
where $(v_k,\Chi_k,\cp_k,\lambda_k):=
(v^{h_k},\Chi^{h_k},\cp^{h_k},\lambda^{h_k})$, $k\in\N$.  

\end{proposition}
\subsection{Convergence in \eqref{eq:td-ns} and \eqref{eq:td-mu}}
We first verify the equations in the bulk.
\begin{proposition}
Let $v,\cp,\Chi$ be the limits obtained in Proposition
\ref{prop:comp-td}. 
Then 
\begin{eqnarray}\nonumber
  \lefteqn{\int_{\Omega_T}\left(-v\partial_t \varphi  + v\cdot\nabla v \varphi
  +\nu(\Chi)Dv:D\varphi \right) \,d(x,t)}\\
&& -
  \int_\Omega \varphi(0,x)\cdot v_0(x)\,dx 
   =\, -\int_{\Omega_T}\Chi\nabla\cp\cdot \varphi\,d(x,t), \label{eq:ns3}
\end{eqnarray}
holds for all $\varphi \in C^\infty([0,T]; C_{0,\sigma}^\infty(\Omega)) $ with
$\varphi|_{t=T}=0$ and
\begin{align}
  &\int_{\Omega_T} \nabla\cp\cdot\nabla\zeta \,d(x,t) \,=\,
  \int_{\Omega_T}\partial_t\zeta\Chi + \dive(\zeta v)\Chi \,d(x,t) +
  \int_\Omega \zeta(0,x)\Chi_0(x)\,dx \label{eq:mu3} 
\end{align}
holds for all $\zeta\in
C^\infty([0,T]\times \overline{\Omega})$ with $\zeta|_{t=T}=0$. 
\end{proposition}
\begin{proof}
If we test in \eqref{eq:td-ns} with $\varphi(\cdot,t)$, where 
$\varphi\in C^\infty([0,T]; C_{0,\sigma}^\infty(\Omega))$ with $\varphi|_{t=T}=0$, and
integrate over $t\in (0,T)$, we  deduce that  
\begin{eqnarray*}
  \lefteqn{\int_{\Omega_T} \Big(-v^h \partial_{t,h}^+\varphi  +
  v^h(.-h)\cdot \nabla v^h\varphi+ \nu(\Chi^h)Dv^h:D\varphi\Big)\sd (x,t)}\\
&&  -\frac{1}{h}\int_0^h\int_\Omega \varphi(x,t) v_0(x)\,dx\,dt
  =\, - \int_{\Omega_T} \Chi^h\nabla\cp^h\cdot \varphi \sd (x,t)
\end{eqnarray*}
for all sufficiently small $h>0$, where we set $\varphi(t)=0$ for $t> T$.
By \eqref{eq:conv-v-strong}, \eqref{eq:conv-v-weak}, \eqref{eq:conv-mu}
we can pass to the limit $h\to 0$ in this equality and obtain 
\begin{eqnarray}\nonumber
  \lefteqn{\int_{\Omega_T} \Big(-v^h\partial_t\varphi + v^h\cdot \nabla v^h
    +\nu(\Chi)Dv:D\varphi\sd (x,t)\varphi\Big) \sd (x,t)}\\ 
&&  -\int_\Omega
  \varphi(0,x) v_0(x) \,dx \,=\, -\int_{\Omega_T} \Chi \nabla \cp \cdot
  \varphi\sd (x,t). \label{eq:td-weak-ns} 
\end{eqnarray}
Similarly we obtain from \eqref{eq:td-mu} that for all $\zeta\in
C^\infty([0,T]\times \ol{\Omega})$ with $\zeta|_{t=T}=0$ 
\begin{eqnarray*}
        \lefteqn{\int_{\Omega_T} \nabla\cp^h\cdot \nabla \zeta \sd (x,t)}\\
         &= &\int_{\Omega_T}\left(\Chi^h\partial_{t,h}^+\zeta +
         (v^h\Chi^h)(.-h)\cdot \nabla \zeta\right) \,d(x,t) 
        +\int_{-h}^0\int_\Omega \Chi_0(x)\zeta(x,t)\,dx\,dt 
\end{eqnarray*}
holds and again we can pass to the limit in this equality and obtain
\eqref{eq:mu3}. 
\end{proof}

\subsection{Convergence in the Gibbs--Thomson law}
The main difficulty in passing the approximate solutions to a limit is
the convergence in the Gibbs--Thomson condition \eqref{eq:td-gt}. In
particular, we cannot exclude that parts of the phase boundary  
$\partial^*\{\Chi^h(\cdot,t)=1\}$ cancel in the limit $h\to 0$. To
overcome such difficulties we consider the limit of the phase boundaries
in the sense of measures and use varifold theory. For the definition of
varifolds and weak mean curvature for varifolds we refer to \cite{Simo83}.

Let $\vf^h_t := \,d|\nabla\Chi^h(\cdot,t)|$ denote the surface
measure of the phase interface $\partial^*\{\Chi^h(\cdot,t)=1\}$,
\begin{gather}
        \vf^h_t(\eta) := \int_\Omega
        \eta\,\,d|\nabla\Chi^h(\cdot,t)|\qquad\text{ for }\eta\in
        C_0(\Omega) \label{eq:def-vf}  
\end{gather}
and let $n^h(t)$ denote the inner normal of
$\partial^\ast\{\Chi^h(\cdot,t)=1\}$, i.e.,
\begin{gather*}
  n^h(x,t)\,=\,
  \frac{\nabla\Chi^h(\cdot,t)}{|\nabla\Chi^h(\cdot,t)|}(x),
\end{gather*}
which is well-defined for $\Ha^{d-1}$-almost all $x\in\partial^*\{\Chi^h(\cdot,t)=1\}$. 
By \eqref{eq:lem-td-gt} the first variation of $\vf^h_t$ is given as
\begin{eqnarray*}
  \delta\vf^h_t(\eta)&=&\int_\Omega
  \big(\dive\eta-n^h(\cdot,t)\cdot
  D\eta\,n^h(\cdot,t)\big)\,d|\nabla\Chi^h_t|\\ 
  &=& \int_\Omega \Chi^h(t)
  \,\dive\big((\cp^h_0(\cdot,t)+\lambda^h(t))\eta\big) \,dx
\end{eqnarray*}
for all $\eta\in C^1_0(\Omega,\Rn)$.

We will prove that for almost all $t\in (0,T)$ the phase boundary
$\partial^*\{\Chi(\cdot,t)=1\}$ has a \emph{generalized mean curvature} in the
following sense.
\begin{definition}\label{def:gen-mc}
Let $E\subset\Omega$ and $\Chi_E\in \BV(\Omega)$. If there exists an
integral\/ $(d-1)$-varifold $\vf$ on\/ $\Omega$ such that
\begin{gather*}
  \partial^*E \,\subset\, \spt (\vf),\\
  \vf\text{ has weak mean curvature vector }{H}_{\vf},\\
  {H}_{\vf}\,\in\, \Lp^s_{\loc}(\vf),\,s>d-1,\,s\geq 2
\end{gather*}
then we call
\[
{H} := {H}_{\vf}|_{\partial^*E}
\]
the generalized mean curvature vector of $\partial^*E$.
\end{definition}
This definition was justified in \cite{Roeg04}, where it is shown that
under the above conditions $H$ is a property of $E$ and independent of
the choice of $\vf$. Moreover, for any $C^2$-hypersurface
$M\subset\Rn$ the mean curvature $H$ of $\partial^*E$ coincides
$\Ha^{d-1}$-almost everywhere on $M\cap \partial^*E$ with the mean
curvature of $M$.
\begin{lemma}\label{lem:conv-vf}
Let $s=4$ if $d=3$ and $1\leq s<\infty$
arbitrary if $d=2$. Then for almost all $t\in (0,T)$ the phase boundary
$\partial^*\{\Chi(\cdot,t)=1\}$ has a generalized mean curvature $H(t)\in
L^s(\,d|\nabla\Chi^h(\cdot,t),\Rn)|$,
\begin{gather}\label{eq:est-mc}
        \int_\Omega |{H}(\cdot,t)|^s \,d|\nabla\Chi(\cdot,t)| \,\leq\,
        C\liminf_{h\to 0}
        \|\cp^{h}(\cdot,t)\|_{H^1(\Omega)}.
\end{gather}
Further, $H(\cdot,t)$ determines the limit of the first variations
$\delta\vf_t^h$:  For any subsequence $h_i\to 0 \;(i\to\infty)$ of $h\to
0$ such that  
\begin{gather}
        \limsup_{i\in\N}
        \|\cp^{h_i}(\cdot,t)\|_{H^1(\Omega)}
        \,<\, \infty \label{eq:subseq-cp} 
\end{gather}
and for all $\eta\in C^1_0(\Omega,\Rn)$ we obtain
\begin{gather}
        \delta\vf^{h_i}_t(\eta)\,\to\,\int_\Omega -{H}(t)\cdot \xi
        \,\,d|\nabla\Chi(\cdot,t)| \qquad\text{ as
        }i\to\infty. \label{eq:conv-Tt} 
\end{gather}
\end{lemma}
\begin{proof}
By Fatou's Lemma and \eqref{eq:est-td} we deduce that
$t\mapsto \liminf_{h\to 0}\|\cp^h(\cdot,t)\|_{H^1(\Omega)}$
belongs 
to $L^2(0,T)$ and that the right-hand side of \eqref{eq:est-mc} is
finite for almost all $t\in (0,T)$. In the following let $t\in (0,T)$ be such
that $\liminf_{h\to 0}\|\cp^h(\cdot,t)\|_{H^1(\Omega)}$ is finite.

Since $\per{\Chi^h(\cdot,t)}$ is uniformly bounded by
\eqref{eq:est-td} and recalling \eqref{eq:subseq-cp} we can extract a
subsequence (not relabeled) $h_i\to 0 (i\to\infty)$ such that 
\begin{align}
  \vf^{h_i}_t \,&\to\, \vf_t \quad\text{ as Radon measures,}
  \label{eq:subsec-vf}\\ 
  \cp^{h_i}(\cdot,t) \,&\schwto\, w_t
  \quad\text{ in } H^1(\Omega) \label{eq:subsec-cp}, 
\end{align}
for a Radon measure $\vf_t$ on $\Omega$ and $w_t\in H^1(\Omega)$. 
We then deduce from \cite[Theorem~1.1]{Scha01} that 
\begin{gather*}
  \vf_t \quad\text{ is an integral varifold},\\
  \vf_t^{h_i}\,\to\, \vf_t\text{ as varifolds, for a subsequence }h_i\to
  0\, (i\to \infty),\\
  \vf_t \quad\text{ has weak mean curvature } {H}_{\vf_t}\in
  L^s(\vf_t), \\
  |\nabla \Chi(\cdot,t)|\,\leq\, \vf_t,
\end{gather*}
and that
\begin{gather*}
  {H}_{\vf_t} = w_tn(t)
\end{gather*}
holds $\vf_t$-almost everywhere, with
\begin{gather*}
  n(\cdot,t) =
  \begin{cases}
    \frac{\nabla\Chi(\cdot,t)}{|\nabla\Chi(\cdot,t)|} & \text{ on }
    \partial^*\{\Chi(\cdot,t)=1\},\\ 
    0 & \text{ elsewhere.}
  \end{cases}
\end{gather*}
According to \cite{Roeg04}
${H}(\cdot,t):={H}_{\vf_t}|_{\partial^*\{\Chi(\cdot,t)=1\}}$ is a
property of $\Chi(\cdot,t)$ and independent of the 
choice of subsequence in \eqref{eq:subsec-vf}, \eqref{eq:subsec-cp}. Moreover,
due to \cite[Theorem~1.2]{Scha01}, we have 
\begin{gather*}
  {H}(\cdot,t) = 0\quad\vf_t\text{-almost everywhere in }
  \{\theta^{d-1}(\vf_t,\cdot)\neq 1\},  
\end{gather*}
where $\theta^{d-1}(\vf_t,\cdot)$ denotes the $(d-1)$-dimensional density of
$\vf_t$, cf. \cite{Simo83}.
Since the first variation is continuous under varifold-convergence, we
then obtain that 
\begin{gather*}
  \lim_{i\to\infty}\,\delta\vf^{{h}_i}_t(\eta) \,=\, \delta\vf_t(\eta)
  =\, \int_\Omega - {H}_{\vf_t}\cdot\eta\,d\vf_t
  =\, \int_\Omega  -{H}(\cdot,t)\cdot\eta \,\,d|\nabla\Chi(\cdot,t)|.
\end{gather*}
\end{proof}
We still have to relate the generalized mean curvature $H(\cdot,t)$ that
we obtained for almost all $t\in (0,T)$ with the weak limit $\cp$
of $\cp^h$ in $L^2(0,T;H^1(\Omega))$.
\begin{lemma}
For all $\xi\in \Lp^2(0,T;C^1_0(\Omega;\Rn))$
\begin{gather}
  \int_0^T -H(\cdot,t)\cdot \xi(\cdot,t) \,d|\nabla\Chi(\cdot,t)|\,dt
  = \int_{\Omega_T}
  \Chi(x,t)\,\dive\big(\cp(x,t)\xi(x,t)\big)\,d(x,t).  \label{eq:gt-int}
\end{gather}
In particular, for almost all $t\in (0,T)$,
\begin{gather}
  H(\cdot,t) = \cp(\cdot,t)n(\cdot,t) \label{eq:gt-final}
\end{gather}
holds $\Ha^{d-1}$-almost everywhere on $\partial^*\{\Chi(\cdot,t)=1\}$.
\end{lemma}
\begin{proof}
Define for all $t\in (0,T)$ such that $H(\cdot,t)\in
L^s(|\nabla\Chi(\cdot,t)|)$ exists 
\begin{gather}
  \dprod{T(t),\psi} \,:=\, \int_\Omega -H(x,t)\cdot
  \psi \,d|\nabla\Chi(x,t)| \label{eq:def-Tt}
\end{gather}
for all $\psi\in C^1_0(\Omega,\Rn)$. Then
\begin{gather}
  |\dprod{T(t),\psi}| \,\leq\,
  C\|\psi\|_{C_0(\Omega)}\|H(\cdot,t)\|_{L^s(\,d|\nabla\Chi^h_t|)}
  \label{eq:est-Tt} 
\end{gather}
and we deduce from \eqref{eq:est-mc} that $T\in L^2(0,T;C_0(\Omega)^*)$.

Similarly we define for $\psi\in C^1_0(\Omega)$
\begin{gather*}
  \dprod{T^h(t),\psi} \,:=\, \delta\vf^h_t(\psi).
\end{gather*}
From \eqref{eq:td-gt} we deduce that
\begin{align}
  \Big| \dprod{T^h(t),\xi(\cdot,t)}\Big|
  \,&=\, \Big|\int_\Omega \Chi^h\nabla \cdot(\xi\cp^h) \,dx \Big| \notag\\
  &\leq\,
  C\|\cp^h(t)\|_{H^1(\Omega)}\|\xi(t)\|_{C^1_0(\Omega)}, 
  \label{eq:est-Tht1}
\end{align}
and
\begin{gather}
  \|T^h\|_{L^2(0,T;C^1_0(\Omega)^\ast)}\,\leq\,
  C\|\cp^h\|_{L^2(0,T;H^1(\Omega))}\|\xi\|_{L^2(0,T;C^1_0(\Omega))}.
  \label{eq:est-Tht}
\end{gather}
Moreover, by \eqref{eq:conv-mu}, \eqref{eq:conv-lambda} there exists a
subsequence $h\to 0$ such that
\begin{gather}
  \lim_{h\to 0} \int_0^T \dprod{T^h(t),\xi(\cdot,t)} \,dt
  \,=\, \int_0^T\int_\Omega \Chi\nabla
  \cdot(\xi\cp) \,dx\,dt. \label{eq:conv-Tht}
\end{gather}
For $\alpha>0$ we now define functions $T^h_\alpha: (0,T)\to
C^1_0(\Omega;\Rn)^*$:
\begin{gather}\label{eq:def-Th}
  \dprod{T^h_\alpha(t),\psi} :=
  \begin{cases}
    \dprod{T^h_t,\psi} &\text{ if }\|\cp^h(\cdot,t)\|_{H^1(\Omega)}\,\leq\,\alpha,\\ 
    \dprod{T(t),\psi} &\text{ else.}       
  \end{cases}
\end{gather}
Fix an arbitrary $\xi\in \Lp^2(0,T;C^1_0(\Omega;\Rn))$. Then we deduce
from Lemma \ref{lem:conv-vf} that for almost all $t\in (0,T)$
\begin{gather}\label{eq:conv-Talt}
  \dprod{T^h_\alpha (t),\xi(\cdot,t)} \,\to\, \int_\Omega
  -H\cdot\xi(\cdot,t) \,\,d|\nabla\Chi(\cdot,t)|. 
\end{gather}
We also see from \eqref{eq:est-Tht1} that
\begin{gather}
  |\dprod{T^h_\alpha (t),\xi(\cdot,t)}| \,\leq\,
  \|\xi(\cdot,t)\|_{C^1_0(\Omega;\Rn)} C\Bigl(\alpha +
  \|T(t)\|_{C_0(\Omega)^*}\Bigr), \label{eq:dom-Thal}
\end{gather}
which gives by \eqref{eq:est-Tt} an $L^1(0,T)$-dominator for the
left-hand side. By \eqref{eq:conv-Talt}, \eqref{eq:dom-Thal} and Lebesgues
Dominated Convergence Theorem we deduce that for all $\alpha>0$
\begin{gather}\label{eq:conv-Th-al}
  \int_0^T \dprod{T^h_\alpha (t),\xi(\cdot,t)}\,dt \,\to\, \int_0^T
  \dprod{T(t)\xi(\cdot,t)}\,dt \quad \text{ as }h\to 0.
\end{gather}
Next, consider the sets $A^h :=\{ t\in (0,T)\,:\,
\|\cp^h(\cdot,t)\|_{H^1(\Omega)}>\alpha\}$ and observe 
\begin{eqnarray*}
  \lefteqn{\Big|\int_0^T \dprod{T^h(t)-T^h_\alpha(t),\xi(\cdot,t)}\,dt \Big| 
  \leq \int_{A^h}\Big|\dprod{T^h(t)-T(t),\xi(\cdot,t)} \Big|\,dt}\\ 
  &\leq& \Bigl(\int_{A^h}\|\xi(\cdot,t)\|_{C^1_0(\Omega;\Rn)}^2 dt
  \Bigr)^{\frac{1}{2}} 
  \Bigl(\|T^h\|_{\Lp^2(0,T;C^1_0(\Omega;\Rn)^*)}+\|T\|_{\Lp^2(0,T;C_0(\Omega;\Rn)^*)}\Bigr). 
\end{eqnarray*}
By \eqref{eq:est-Tt} and \eqref{eq:est-Tht}
$\|T^h\|_{\Lp^2(0,T;C^1_0(\Omega;\Rn)^*)}$ and 
$\|T\|_{\Lp^2(0,T;C^1_0(\Omega;\Rn)^*)}$ are bounded uniformly in $h>0$.
Since
\begin{gather*}
 |A^h|\,\leq\,\frac{1}{\alpha^2}\|\cp^h\|_{\Lp^2(0,T;H^1(\Omega))}^2 
 \,\leq\, \frac{1}{\alpha^2}C,    
\end{gather*}
we end up with
\begin{gather}\label{eq:conv-T-diff}
  \int_0^T \dprod{T^h(t)-T^h_\alpha (t),\xi(\cdot,t)}\,dt \,\to\, 0\qquad
  \text{as}\ \alpha\to\infty
\end{gather}
uniformly in $h>0$. Thus we obtain from \eqref{eq:conv-Talt},
\eqref{eq:conv-T-diff} that
\begin{align}
  \int_{\Omega_T} \Chi\dive(\cp\xi)\,d(x,t) =\,& \lim_{h\to 0}
  \int_{\Omega_T}
  \Chi^h\dive(\cp^h\xi)\,d(x,t)\\ 
  =\,& \lim_{h\to 0}\int_0^T \dprod{\xi(\cdot,t),\delta\vf^h_t}\,dt
  =\, \int_0^T \dprod{\xi(\cdot,t),T(t)}\,dt,
\end{align}
which proves \eqref{eq:gt-int}.
Since no time derivative is involved, we deduce that for almost all $t\in
(0,T)$ and all $\xi\in C^1_0(\Omega;\Rn)$, 
\[
\int_{\Omega} \Chi(\cdot,t)\dive\big(\cp(\cdot,t)\xi\big)\sd x = \dprod{T(t),\xi}
\]
holds. The Gauss--Green theorem \cite[Theorem 5.8.1]{EGar92} and \eqref{eq:def-Tt} yield
\begin{gather*}
  \int_\Omega \cp(\cdot,t)n(\cdot,t)\cdot\xi\,\,d|\nabla\Chi(\cdot,t)| =
  \int_\Omega {H}(\cdot,t)\cdot\xi\, \,d|\nabla\Chi(\cdot,t)|,
\end{gather*}
with $n(\cdot,t)= \frac{\nabla\Chi(\cdot,t)}{|\nabla\Chi(\cdot,t)|}$ on $\partial^*\{\Chi(\cdot,t)=1\}$.
This finally proves \eqref{eq:gt-final}.
\end{proof}

\begin{appendix}
\section{Sharp Interface Limit}

Here we discuss the relation between (\ref{eq:1})-(\ref{eq:6}) and its diffuse
interface analog (\ref{eq:NSCH1})-(\ref{eq:NSCH7}). First we consider
the corresponding energy identities. For the
Navier--Stokes/Mullins--Sekerka system we recall that by
\eqref{eq:en_id_nsms} every sufficiently
smooth solution of  (\ref{eq:1})-(\ref{eq:6}) satisfies 
\begin{equation}\label{eq:en_id_nsms'}
  \frac{d}{dt} \frac12 \int_\Omega |v(t)|^2\sd x +\STC \frac{d}{dt}
  \Ha^{d-1}(\Gamma(t)) = - \int_\Omega \nu(\Chi) |Dv|^2 \sd x - m\int_\Omega |\nabla \mu|^2 \sd x. 
\end{equation}
On the other hand, every sufficiently smooth solution of (\ref{eq:NSCH1})-(\ref{eq:NSCH7}) satisfies
\begin{equation}\label{eq:en_id_nsch}
  \frac{d}{dt} \frac12 \int_\Omega |v(t)|^2\sd x +\frac{d}{dt} E_\eps(c(t)) = - \int_\Omega \bar{\nu}(c) |Dv|^2 \sd x -m\int_\Omega |\nabla \mu|^2 \sd x, 
\end{equation}
where 
\begin{equation*}
  E_\eps(c) = \frac{\eps}2\int_\Omega |\nabla c|^2 \sd x + \frac1\eps \int_\Omega f(c) \sd x.
\end{equation*}
Moreover, by Modica and Mortola~\cite{ModicaMortola1} or Modica~\cite{Modica1}, 
we have
\begin{equation*}
  E_\eps \to_{\eps\to 0} \mathcal{P} \qquad \text{w.r.t.}\ L^1\text{-}\Gamma\text{-convergence},
\end{equation*}
where
\begin{equation*}
  \mathcal{P}(u)=
  \begin{cases}
    \STC \Ha^{d-1}(\partial^\ast E) & \text{if}\ u = 2\Chi_E -1 \ \text{and}\ E\ \text{has finite perimeter},\\
+\infty& \text{else}.
  \end{cases}
\end{equation*}
Here $\STC= \int_{0}^1\sqrt{2f(s)} \sd s$, $f\colon \R\to [0,\infty)$ is a suitable function such that $f(s)=0$ if and only if $s=0,1$, 
and $\partial^\ast E$ denotes the reduced boundary. Note that $\partial^\ast E= \partial E$ if $E$ is a sufficiently regular domain. 
Therefore we see that for constant  $m>0$  the energy identity
\eqref{eq:en_id_nsms'} is formally identical to the sharp interface
limit of the energy identity \eqref{eq:en_id_nsch} of the diffuse interface
model (\ref{eq:NSCH1})-(\ref{eq:NSCH7}).
In contrast, if we would choose $m=m_\eps\to_{\eps\to 0} 0$ in
(\ref{eq:NSCH1})-(\ref{eq:NSCH7}), we formally obtain in
the sharp interface limit the energy
identity of the classical two-phase flow (\ref{eq:1})-(\ref{eq:6}). 

Now we will adapt the arguments of Chen~\cite{Chen96} to show that as
$\eps\to 0$ and if $m=m(\eps)\to_{\eps\to 0} m_0 >0$
solutions of the diffuse interface model (\ref{eq:NSCH1})-(\ref{eq:NSCH7})
converge to \emph{varifold solutions} of the system (\ref{eq:1})-(\ref{eq:6}),
which are defined as follows:
\begin{definition}\label{defn:VarifoldSolution}
 Let $v_0\in L^2_\sigma(\Omega)$, $E_0\subset \Omega$ be a set of finite
 perimeter, and let $Q=\Omega\times (0,\infty)$.  Then $(v,E,\mu,V)$ is called \emph{varifold solution} of
  (\ref{eq:1})-(\ref{eq:6}) with initial values $(v_0,E_0)$ if the following conditions are satisfied:
  \begin{enumerate}
  \item $v\in L^\infty(0,\infty;L^2_\sigma(\Omega))\cap
  L^2(0,\infty;H^1_0(\Omega)^d)$, $\mu\in L^2_{loc}([0,\infty);H^1(\Omega))$
  and $\nabla \mu \in L^2(Q)$.
\item $E=\bigcup_{t\geq 0} E_t\times \{t\}$ is a measurable subset of
  $\Omega\times [0,\infty)$ such that $\Chi_E \in
  C([0,\infty);L^1(\Omega))\cap L^\infty(0,\infty;BV(\Omega))$, $\int_\Omega
  \Chi_E(t) \sd x = m_0$ for all $t\geq 0$, and
  $\Chi_E|_{t=0}=\Chi_{E_0}$ in $L^1(\Omega)$.
\item $V$ is a Radon measure on $\Omega\times G_{d-1}\times (0,\infty)$ 
  such that $V= V^t dt$ where $V^t$ is for almost all $t\in (0,\infty)$ a
  general varifold on $\Omega$, i.e., a Radon measure on $\Omega \times
  G_{d-1}$. Moreover, for almost all $t\in (0,\infty)$ $V^t$ has the representation
  \begin{equation*}
    \int_{\Omega \times G_{d-1}} \psi (x,p) \sd V^t(x,p)  = \sum_{i=1}^d \int_\Omega
    b_i^t(x) \psi(x,p_i^t(x))\sd \lambda^t(x) 
  \end{equation*}
  for all $\psi\in
    C_0(\Omega\times G_{d-1})$,
  where $\lambda^t$ is a Radon measure on $\overline{\Omega}$, $b_1^t,\ldots,
  b_N^t$ and $p_1^t,\ldots, p_N^t$ are measurable $\R$- and $\Rn$-valued
  functions, respectively, such that 
  \begin{alignat*}{1}
    &0\leq b_i^t \leq 1,\quad \sum_{i=1}^d b_i^t \geq 1,\quad \sum_{i=1}^d
    p_i^t\otimes p_i^t = I \quad \mu^t\text{-almost everywhere}\\
    &\frac{|\nabla \Chi_{E_t}|}{\lambda^t} \leq \frac1{\STC}\qquad
    \lambda^t\text{-a.e. with}\ \STC= \int_{0}^1 \sqrt{2f(s)} \sd s.
  \end{alignat*}
\item Moreover,  
\begin{eqnarray}\nonumber
  \lefteqn{\int_{Q}\left(-v\partial_t \varphi  +  v\cdot\nabla v \varphi
  +\nu(\Chi_E)Dv:D\varphi\right)  \,d(x,t)}\\
&& -
  \int_\Omega \varphi|_{t=0}\cdot v_0\,dx 
   = -\STC \int_Q \Chi_E\nabla \mu\cdot \varphi \sd (x,t) \label{eq:ns'}
\end{eqnarray}
holds for all $\varphi \in C^\infty([0,\infty); C_{0,\sigma}^\infty(\Omega)) $
with $\supp \varphi \subset \Omega\times [0,T]$ for some $T>0$,
\begin{align}
  &m\int_{Q} \nabla\cp\cdot\nabla\psi \,d(x,t) \,=\,
  \int_{Q}\Chi(\partial_t\psi + \dive(\psi v)) \,d(x,t) +
  \int_\Omega \psi|_{t=0}\Chi_{E_0}\,dx \label{eq:mu'} 
\end{align}
holds for all $\psi\in
C^\infty([0,\infty)\times \overline{\Omega})$ with $\supp\psi \subset\ol{\Omega}\times [0,T]$
for some $T>0$, as well as
\begin{equation*}
  (\Chi_{E_t}, \Div (\mu\eta))_\Omega =  \left\langle \delta
    V^t,\eta\right\rangle :=
  \int_{\Omega\times G_{d-1}} (I-p\otimes p):\nabla \eta (x) \sd V^t(x,p) 
\end{equation*}
for all $\eta \in C^1_0(\Omega \times G_{d-1})$.
\item Finally, for almost all $0\leq s\leq t <\infty$
  \begin{eqnarray}\nonumber
    \lefteqn{\frac12\|v(t)\|_{L^2(\Omega)}^2+\lambda^t (\ol{\Omega})}\\\label{eq:EnergyEstim}
&& +
    \int_s^t\int_\Omega \left(\nu(\Chi_{E_\tau}) |Dv|^2 + m |\nabla\mu|^2
    \right)\sd x\sd \tau \leq \frac12\|v(s)\|_{L^2(\Omega)}^2+\lambda^s (\ol{\Omega}).
  \end{eqnarray}
\end{enumerate}
\end{definition}
Here and in the following $f$ shall satisfy $f\in C^3(\R)$, $f(s)\geq 0$ and $f(s)=0$ if and only if
$s=0,1$ as well as $f''(s)\geq c_0 (1+|s|)^{p-2}$ if $s\geq 1-c_0$ and if
$s\leq c_0$ for some
$c_0>0, p>2$. Then $F(s):= f(\frac{s+1}2)$ we will satisfy the assumption in \cite{Chen96}. We will even assume that $p\geq 3$, which we will need in the
following to estimate $v_\eps \cdot \nabla c_\eps = \Div (v_\eps c_\eps)$
uniformly in $L^2(0,T;H^{-1}(\Omega)$.  One can choose e.g. $f(s) = s^2(1-s)^2$.

For the following we denote 
\begin{equation*}
  e_\eps (c) = \frac{\eps}2 |\nabla c|^2 + \eps^{-1} f(c).
\end{equation*}

\begin{theorem}\label{thm:SharpInterfaceLimit}
  Let $\Omega\subset\R^d$, $d=2,3$, be a smooth bounded domain and let $\nu\in
  C^0(\R)$ with $\nu(s)\geq \nu_0 >0$ for all $s\in\R$.  Moreover, let
  initial data $(v_{0,\eps},c_{0,\eps})\in
  L^2_\sigma(\Omega)\times H^1(\Omega)$ with $\frac1{|\Omega|}\int_\Omega
  c_{0,\eps} \sd x=
  \bar{c} \in (-1,1)$  be given that satisfy
  \begin{equation}\label{eq:EstimIV}
    \frac12 \int_\Omega |v_{0,\eps}|^2 \sd x + E_\eps(c_{0,\eps}) \leq R 
  \end{equation}
  uniformly in $\eps\in (0,1]$ for some $R>0$. Finally, let
$(m_\eps)_{\eps\in (0,1]}$ with $m_\eps\to_{\eps\to 0} m >0$. Consider
now (\ref{eq:NSCH1})-(\ref{eq:NSCH7}) with $m$ 
  replaced by $m_\eps$ and let for every $\eps\in (0,1]$ 
  $(v_\eps,c_\eps,\mu_\eps)$ be weak solutions in the sense of
  \cite[Definition~1]{ModelH}. 
  Then there is a subsequence
  $(\eps_k)_{k\in\N}$ such that $\eps_k \searrow 0$ as $k\to \infty$ such that:
  \begin{enumerate}
  \item There are measurable sets $E\subset \Omega \times [0,\infty)$ and
    $E_0\subset \Omega$ such that
    $c_{\eps_k}\to_{k\to\infty} \Chi_{E}$ almost everywhere in $\Omega \times
  [0,\infty)$ and in $C^{\frac19}([0,T];L^2(\Omega))$ for all $T>0$ as well as
  $c_{0,\eps_k}\to_{k\to\infty} -1+2\Chi_{E_0}$ almost everywhere in $\Omega$
  and $\Chi_{E_0}= \Chi_{E}|_{t=0}$ in $L^2(\Omega)$.
\item There are $\mu\in L^2_{loc}([0,\infty);H^1(\Omega)),
  v\in L^2(0,\infty;H^1_0(\Omega))$, $v_0\in L^2_\sigma(\Omega)$ such that
  \begin{alignat}{2}\label{eq:Conv1}
    \mu_{\eps_k} &\rightharpoonup_{k\to \infty} \mu &\qquad& \text{in}\
    L^2(0,T;H^1(\Omega))\quad \text{for all}\ T>0,\\\label{eq:Conv2}
    v_{\eps_k} &\rightharpoonup_{k\to \infty} v &\qquad& \text{in}\
    L^2(0,\infty;H^1(\Omega)),\\\label{eq:Conv3}
   v_{0,\eps_k} &\rightharpoonup_{k\to \infty} v_0 &\qquad& \text{in}\
    L^2_\sigma (\Omega).
  \end{alignat} 
\item There exist Radon measures $\RM$ and $\RM_{ij}$, $i,j=1,\ldots, n$ on
  $\ol{\Omega}\times [0,\infty)$ such that for every $T>0$, $i,j=1,\ldots, n$
  \begin{alignat}{2}\label{eq:RM1}
    e_{\eps_k} (c_{\eps_k}) \sd x\sd t  &\rightharpoonup_{k\to \infty}^\ast 
    \RM &\qquad& \text{in}\ \mathcal{M}(\ol{\Omega}\times [0,T]),\\\label{eq:RM2}
    \eps_k\partial_{x_i}c_{\eps_k}\partial_{x_j} c_{\eps_k} \sd x\sd t &\rightharpoonup_{k\to \infty}^\ast 
    \RM_{ij} &\qquad& \text{in}\ \mathcal{M}(\ol{\Omega}\times [0,T]).
  \end{alignat}
\item There exist a Radon $V= V^t \sd t$ on $\Omega\times G_{d-1}\times
  (0,\infty)$ such that $(v,E,\mu,V)$ is a varifold solution of
  (\ref{eq:1})-(\ref{eq:6}) with initial values $(v_0,E_0)$ in the sense of
  Definition~\ref{defn:VarifoldSolution}, where
  \begin{equation}\label{eq:ReprFirstVar}
    \int_0^T \weight{\delta V^t,\eta} \sd t = \int_0^T \int_\Omega
    \nabla \eta : \left(\sd \RM I- (d\RM_{ij})_{i,j=1}^d\right)\sd t
  \end{equation}
  for all $\eta \in C^1_0(\Omega \times [0,T];\R^d)$. 
  \end{enumerate}
\end{theorem}
First of all, by the definition of weak solutions in \cite{ModelH}, we have
  \begin{eqnarray}\nonumber
    \lefteqn{\int_Q \left(-v_\eps\partial_t \psi+ v_\eps\cdot \nabla v_\eps\psi + \nu(c_\eps)Dv_\eps: D \psi\right) \sd
    (x,t)}\\\label{eq:weakNSCH1}
&&- \int_\Omega v_0\psi|_{t=0}\sd x  = -\int_Q 
    c_\eps\nabla \mu_\eps  \cdot \psi\sd (x,t)
  \end{eqnarray}
  for all $\psi \in C^\infty([0,\infty)\times\Omega)^d$ with $\Div \psi=0$ and
  $\psi(t) =0$ for $t\geq T$ for some $T>0$, as well as
  \begin{alignat}{1}\label{eq:weakNSCH2}
    m_\eps \int_Q\nabla\mu_\eps \cdot \nabla \varphi \sd(x,t) &= \int_Q
    c_\eps \left(\partial_t \varphi +\Div (\varphi v_\eps)\right)\sd (x,t) +
    \int_\Omega c_{0,\eps}\varphi|_{t=0}\sd x \\\label{eq:weakNSCH3}
    \int_Q \mu_\eps\varphi \sd (x,t) &= \int_Q f(c_\eps)\varphi\sd (x,t) + \int_Q
    \nabla c_\eps\cdot \nabla \varphi\sd (x,t)
  \end{alignat}
  for all $\varphi \in C^\infty([0,\infty)\times\ol{\Omega})$ with $\supp
  \varphi\subset \ol{\Omega}\times[0,T]$ for some $T>0$. Moreover, from the
  energy inequality in \cite[Definition~1]{ModelH} we obtain
\begin{equation*}
  \frac12 \|v_\eps(t)\|_{L^2(\Omega)}^2 +E_\eps(c_{\eps}(t))
  +\int_0^t\int_\Omega \left(\nu(c_\eps) |Dv_\eps|^2 +m_\eps |\nabla
    \mu_\eps|^2\right) \sd x \sd t \leq R 
\end{equation*}
for all $t\in [0,\infty)$. Therefore we have
\begin{equation}
  \label{eq:H1ChemPotEstim}
  \|\mu_\eps (\cdot,t)\|_{H^1(\Omega)} \leq C \left(E_\eps(c_\eps(t)) + \|\nabla \mu_\eps(.,t)\|_{L^2(\Omega)}\right)
\end{equation}
for all $t>0$ and $0<\eps\leq \eps_0$ for some $C,\eps_0>0$ due to \cite[Lemma~3.4]{Chen96}.

Hence there exists a subsequence
$\eps_k\searrow 0$ as $k\to\infty$ such that (\ref{eq:Conv1})-(\ref{eq:Conv3})
holds. Moreover, using (\ref{eq:NSCH1}) and the Lemma by Aubin-Lions, one
easily derives that $v_{\eps_k} \to_{k\to \infty} v$ strongly in $L^2(\Omega\times
(0,T))$ for all $T>0$ and $v_{\eps_k}(t) \to_{k\to \infty} v(t)$ strongly in
$L^2(\Omega)$ for almost every $t\geq 0$. 

Using the assumptions on $f$ we further deduce that
\begin{equation*}
  \int_\Omega |c_\eps(t)|^p \sd x \leq C(1+R),\quad \int_\Omega \dist(c,\{0,1\})^2
  \sd x\leq C\eps R 
\end{equation*}
uniformly in $\eps\in (0,1]$, for almost every $0<t<\infty$ and
$t=0$. Now we define as in \cite{Chen96}
\begin{equation*}
  w_\eps (x,t) = W(c_\eps(x,t))\qquad \text{where}\ W(c)=
  \int_{0}^c\sqrt{2\tilde{f}(s)}\sd s,\ \tilde{f}(s)= \min (f(s),1+|s|^2).
\end{equation*}
Then $(\nabla w_\eps)_{\eps\in (0,1]}$ is uniformly bounded in
$L^\infty(0,\infty;L^1(\Omega))$ since
\begin{equation}\label{eq:NablaWepsEstim}
  \int_\Omega |\nabla w_\eps(x,t)|\sd x = \int_\Omega
  \sqrt{2\tilde{f}(s)}|\nabla c_\eps(x,t)|\sd x \leq \int_\Omega e_\eps
  (c_\eps(t))\sd x \leq R.
\end{equation}
Moreover, we have for all $u_1,u_2\in\R$
\begin{equation}\label{eq:WEstim}
  c_1|u_1-u_2|^2 \leq |W(u_1)-W(u_2)|\leq c_2 |u_1-u_2|^2 (1+|u_1|+|u_2|)
\end{equation}
which again follows easily from the assumptions on $f$.
Now we obtain:
\begin{lemma}
  There is some $C$ independent of $\eps\in (0,1]$ such that
  \begin{equation*}
    \|w_\eps\|_{C^{\frac18}([0,\infty); L^1(\Omega))}+
    \|c_\eps\|_{C^{\frac18}([0,\infty); L^2(\Omega))} \leq C.
  \end{equation*}
\end{lemma}
\begin{proof}
  The proof is a modification of \cite[Lemma~3.2]{Chen96}. Therefore we only
  give a brief presentation, describing the differences. 

  For sufficiently small $\eta>0$ let
  \begin{equation*}
    c_\eps^\eta (x,t) = \int_{B_1} \rho(y) c_\eps (x-\eta y,t) \sd y,\qquad
    x\in\Omega, t\geq 0,
  \end{equation*}
  where $\rho$ is a standard mollifying kernel and $c_\eps$ is extended in an
  $\eta_0$-neighborhood of $\Omega$ as in \cite[Proof of Lemma~3.2]{Chen96}.
  Then one obtains
  \begin{alignat*}{1}
    \|\nabla c^\eta_\eps(.,t)\|_{L^2(\Omega)} &\leq
    C\eta^{-1}\|c_\eps(.,t)\|_{L^2(\Omega)}\leq C'\eta^{-1}\\
    \|c^\eta_\eps(.,t) - c_\eps(.,t)\|_{L^2(\Omega)}^2&\leq C\eta \|\nabla
    w_\eps(.,t)\|_{L^1(\Omega)}\leq C'\eta 
  \end{alignat*}
for all sufficiently small $\eta>0$, cf. \cite[Proof of Lemma~3.2]{Chen96}.
Next we use that
\begin{equation*}
  (c_\eps(\cdot,t) - c_\eps(\cdot,\tau),\varphi)_\Omega = -\int_\tau^t \left((\nabla
    \mu- v_\eps
    c_\eps)(s), \nabla \varphi \right)_\Omega \sd s
\end{equation*}
for all $\varphi \in C^1(\ol{\Omega})$
because of (\ref{eq:NSCH3}) in its weak form.  
Here
\begin{equation*}
  \|\nabla\mu- v_\eps c_\eps\|_{L^2(\Omega\times (\tau,t))} \leq C(R)\quad
  \text{for all}\ 0\leq \tau\leq t <\infty, |t-\tau|\leq 1
\end{equation*}
since $v_\eps \in L^2(0,\infty;L^6(\Omega))$ and $c_\eps \in
L^\infty(0,\infty;L^3(\Omega))$ are bounded due to $d\leq 3$ and $p\geq
3$. Hence 
\begin{eqnarray*}
  \lefteqn{\int_\Omega (c_\eps^\eta(x,t)- c_\eps^\eta (x,\tau))(c_\eps(x,t)- c_\eps
  (x,\tau)) \sd x}\\
&=& -\int_\tau^t \left((\nabla
    \mu- v_\eps
    c_\eps)(s), \nabla c_\eps^\eta(x,t)- \nabla c_\eps^\eta (x,\tau) \right)_\Omega
  \sd s\\
&\leq& C(R) (t-\tau)^{\frac12} \sup_{s\in [0,\infty}\|\nabla c_\eps^\eta
(.,s)\|_{L^2(\Omega)} \leq C(T,R) \eta^{-1}(t-\tau)^{\frac12}
\end{eqnarray*}
for all $0\leq \tau\leq t <\infty, |t-\tau|\leq 1$.
Now, choosing $\eta= \min(\eta_0,(t-\tau)^{\frac14})$, we conclude
\begin{equation*}
  \|c_\eps (.,t)-c_\eps(.,\tau)\|_{L^2(\Omega)}^2 \leq
  C|t-\tau|^{\frac14}\qquad \text{for all}\ |t-\tau|\leq 1
\end{equation*}
from the previous estimates. Thus $c_\eps \in
C^{\frac18}([0,\infty);L^2(\Omega))$. Using (\ref{eq:WEstim}), one derives
$w_\eps \in C^{\frac18}([0,\infty);L^1(\Omega))$ as in \cite{Chen96}.
\end{proof}
Note that the previous lemma and (\ref{eq:NablaWepsEstim}) implies that 
$(w_\eps)_{0<\eps\leq 1}$ is bounded in $L^\infty(0,\infty;BV(\Omega))$.
\begin{lemma}
  There is a subsequence $\eps_k\searrow 0$ as $k\to\infty$, a measurable
  function $\mathcal{E}(t)$, $t\in (0,\infty)$ and a measurable set $E\subset
  \Omega\times (0,\infty)$ such that 
  \begin{enumerate}
  \item $E_{\eps_k}(c_{\eps_k}(t))\to \mathcal{E}(t)$ for almost all $t\geq
    0$.
  \item $w_{\eps_k}\to \kappa \Chi_E$ almost everywhere in $\Omega\times
    (0,\infty)$ and in $C^{\frac19}([0,T];L^1(\Omega))$ for all $T>0$. 
  \item $c_{\eps_k}\to  \Chi_E$ almost everywhere in $\Omega\times
    (0,\infty)$ and in $C^{\frac19}([0,T];L^2(\Omega))$ for all $T>0$.
  \end{enumerate}
  In addition, $\Chi_E\in L^\infty(0,\infty;BV(\Omega))\cap C^{\frac14}([0,\infty);L^1(\Omega))$ and $E_t:=\{x\in\Omega: (x,t)\in E\}$
  satisfies $|E_t|= |E_0|= \frac{1+\bar{c}}2|\Omega|$ for almost all $t\geq 0$ and
  \begin{equation*}
    |\nabla \Chi_{E_t}|(\Omega) \leq \frac1{\STC} \mathcal{E}(t)\leq
    \frac1{\STC} \mathcal{E}_0. 
  \end{equation*}
\end{lemma}
\begin{proof}
  The main difference to the proof of \cite[Lemma~3.3]{Chen96} is the proof of
  convergence for $E_{\eps_k}(c_{\eps_k}(t))$. To this end one uses that
  \begin{equation*}
    F_\eps(t):= \frac12 \|v_\eps(t)\|_{L^2(\Omega)}^2 + E_{\eps}(c_\eps(t)),\qquad t\geq 0,
  \end{equation*}
  is a sequence of bounded, monotone decreasing functions and
  $v_{\eps_k}(t)\to_{k\to\infty} v(t)$ for almost all $t\geq 0$ in
  $L^2(\Omega)$. The rest of the proof is identical with the proof of \cite[Lemma~3.3]{Chen96}. 
\end{proof}
Using the previous statements one can now easily finish the proof of
Theorem~\ref{thm:SharpInterfaceLimit} by the arguments of
\cite[Section~3.5]{Chen96}. In particular, (\ref{eq:ns'}) and (\ref{eq:mu'})
easily follow from (\ref{eq:weakNSCH1}) and  (\ref{eq:weakNSCH2}).  It mainly remains to show
(\ref{eq:ReprFirstVar}) and (\ref{eq:EnergyEstim}). Let $\RM, \RM_{i,j}$ be as
in (\ref{eq:RM1})-(\ref{eq:RM2}). To show the energy
estimate (\ref{eq:EnergyEstim}) one uses that 
$
  \sd \RM = \sd \RM^t \sd t  
$ for some Radon measures $\RM^t$ on $\ol{\Omega}$ and that 
for almost every $0\leq t\leq
s<\infty$
\begin{eqnarray*}
  \lefteqn{\RM^t(\ol{\Omega}) = \lim_{k\to\infty} E_{\eps_k}(c_{\eps_k}(t))} \\
  &\leq &\lim_{k\to\infty} \left(E_{\eps_k}(c_{\eps_k(s)}) - \int_s^t
    \int_\Omega \left(\nu(c_{\eps_k}) |D v_{\eps_k}|^2 + m_{\eps_k}
      |\nabla \mu_{\eps_k}|^2\right) \sd x\sd \tau\right)\\
  && + \lim_{k\to \infty} \left(\frac12\|v_{\eps_k}(s)\|_{L^2(\Omega)}^2
    -\frac12\|v_{\eps_k}(t)\|_{L^2(\Omega)}^2 \right)\\
  &\leq & \RM^t(\ol{\Omega}) - \int_s^t
  \int_\Omega \left(\nu(\Chi_E) |D v|^2 + m
    |\nabla \mu|^2 \right)\sd x\sd \tau + \frac12\|v(s)\|_{L^2(\Omega)}^2
    -\frac12\|v(t)\|_{L^2(\Omega)}^2, 
\end{eqnarray*}
where we have used the weak convergence of $Dv_{\eps_k}, \nabla\mu_{\eps_k}$
in $L^2(Q)$ and the strong convergence of $v_{\eps_k}(t)$ in $L^2(\Omega)$ for almost
all $t\geq 0$. Finally, (\ref{eq:ReprFirstVar}) follows in precisely the same
way as in \cite[Section~3.5]{Chen96}, where we note that the main point in the
argumentation is that the discrepancy measure 
$$
\xi^\eps(c_\eps):= \frac{\eps}2|\nabla c_\eps|^2 -\frac1\eps f(c_\eps)
$$
converges to a non-positive measure. The latter fact follows from
\cite[Section~3.5]{Chen96} since we have the same bounds on
$\mu_\eps$ in $H^1(\Omega)$.
\label{app:}

\end{appendix}

\def\cprime{$'$} \def\cprime{$'$}


\begin{thebibliography}{10}

\bibitem{ModelH}
H.~Abels.
\newblock On a diffuse interface model for two-phase flows of viscous,
  incompressible fluids with matched densities.
\newblock {\em Arch. Rat. Mech. Anal., DOI 10.1007/s00205-008-0160-2}.

\bibitem{GeneralTwoPhaseFlow}
H.~Abels.
\newblock On generalized solutions of two-phase flows for viscous
  incompressible fluids.
\newblock {\em Interfaces Free Bound.}, 9:31--65, 2007.

\bibitem{ReviewGeneralTwoPhaseFlow}
H.~Abels.
\newblock On the notion of generalized solutions of two-phase flows for viscous
  incompressible fluids.
\newblock {\em RIMS K{\^{o}}ky{\^{u}}roku Bessatsu}, B1:1--15, 2007.

\bibitem{AmbrosioEtAl}
L.~Ambrosio, N.~Fusco, and D.~Pallara.
\newblock {\em {Functions of bounded variation and free discontinuity
  problems.}}
\newblock {Oxford Mathematical Monographs. Oxford: Clarendon Press. xviii, 434
  p.}, 2000.

\bibitem{Chen96}
X.~Chen.
\newblock Global asymptotic limit of solutions of the {C}ahn-{H}illiard
  equation.
\newblock {\em J. Differential Geom.}, 44(2):262--311, 1996.

\bibitem{DenisovaTwoPhase}
I.~V. Denisova and V.~A. Solonnikov.
\newblock Solvability in {H}\"older spaces of a model initial-boundary value
  problem generated by a problem on the motion of two fluids.
\newblock {\em Zap. Nauchn. Sem. Leningrad. Otdel. Mat. Inst. Steklov. (LOMI)},
  188(Kraev. Zadachi Mat. Fiz. i Smezh. Voprosy Teor. Funktsii. 22):5--44, 186,
  1991.

\bibitem{EGar92}
L.~C. Evans and R.~F. Gariepy.
\newblock {\em Measure theory and fine properties of functions}.
\newblock Studies in Advanced Mathematics. CRC Press, Boca Raton, FL, 1992.

\bibitem{GurtinTwoPhase}
M.~E. Gurtin, D.~Polignone, and J.~Vi{\~n}als.
\newblock Two-phase binary fluids and immiscible fluids described by an order
  parameter.
\newblock {\em Math. Models Methods Appl. Sci.}, 6(6):815--831, 1996.

\bibitem{HohenbergHalperin}
P.~Hohenberg and B.~Halperin.
\newblock Theory of dynamic critical phenomena.
\newblock {\em Rev. Mod. Phys.}, 49:435--479, 1977.

\bibitem{LiuShenModelH}
C.~Liu and J.~Shen.
\newblock A phase field model for the mixture of two incompressible fluids and
  its approximation by a {F}ourier-spectral method.
\newblock {\em Phys. D}, 179(3-4):211--228, 2003.

\bibitem{LStu95}
S.~Luckhaus and T.~Sturzenhecker.
\newblock Implicit time discretization for the mean curvature flow equation.
\newblock {\em Calc. Var. Partial Differential Equations}, 3(2):253--271, 1995.

\bibitem{MaekawaTwoPhaseFlow}
Y.~Maekawa.
\newblock On a free boundary problem for viscous incompressible flows.
\newblock {\em Interfaces Free Bound.}, 9(4):549--589, 2007.

\bibitem{Modica1}
L.~Modica.
\newblock The gradient theory of phase transitions and the minimal interface
  criterion.
\newblock {\em Arch. Rational Mech. Anal.}, 98(2):123--142, 1987.

\bibitem{ModicaMortola1}
L.~Modica and S.~Mortola.
\newblock Un esempio di {$\Gamma \sp{-}$}-convergenza.
\newblock {\em Boll. Un. Mat. Ital. B (5)}, 14(1):285--299, 1977.

\bibitem{PlotnikovTwoPhase}
P.~Plotnikov.
\newblock {Generalized solutions to a free boundary problem of motion of a
  non-Newtonian fluid.}
\newblock {\em Sib. Math. J.}, 34(4):704--716, 1993.

\bibitem{Roeg04}
M.~R{\"o}ger.
\newblock Solutions for the {S}tefan problem with {G}ibbs-{T}homson law by a
  local minimisation.
\newblock {\em Interfaces Free Bound.}, 6(1):105--133, 2004.

\bibitem{Scha01}
R.~Sch{\"a}tzle.
\newblock Hypersurfaces with mean curvature given by an ambient {S}obolev
  function.
\newblock {\em J. Differential Geom.}, 58(3):371--420, 2001.

\bibitem{Simo87}
J.~Simon.
\newblock Compact sets in the space {$L\sp p(0,T;B)$}.
\newblock {\em Ann. Mat. Pura Appl. (4)}, 146:65--96, 1987.

\bibitem{Simo83}
L.~Simon.
\newblock {\em Lectures on geometric measure theory}, volume~3 of {\em
  Proceedings of the Centre for Mathematical Analysis, Australian National
  University}.
\newblock Australian National University Centre for Mathematical Analysis,
  Canberra, 1983.

\bibitem{Sohr01}
H.~Sohr.
\newblock {\em The {N}avier-{S}tokes equations}.
\newblock Birkh\"auser Advanced Texts: Basler Lehrb\"ucher. [Birkh\"auser
  Advanced Texts: Basel Textbooks]. Birkh\"auser Verlag, Basel, 2001.
\newblock An elementary functional analytic approach.

\end{thebibliography}
%
\end{document}